\documentclass[12 pt]{article}
\usepackage[pagebackref]{hyperref}
\hypersetup{
colorlinks = true,
linkcolor = blue,
citecolor = blue}
\usepackage[utf8]{inputenc}
\usepackage[T1]{fontenc}

\usepackage[]{amsfonts}
\usepackage{amssymb}
\usepackage{amsmath,amsthm}
\def\couleur(#1 #2 #3)
	{
\special{ps::[end]}}

\def\bx#1{\setbox1=\hbox{\kern3pt{#1}\kern3pt}			
 \dimen1=\ht1 \advance\dimen1 by 3pt \dimen2=\dp1 \advance\dimen2 by 3pt
 \setbox1=\hbox{\vrule height\dimen1 depth\dimen2\box1\vrule}%
 \setbox1=\vbox{\hrule\box1\hrule}%
 \advance\dimen1 by .4pt \ht1=\dimen1
 \advance\dimen2 by .4pt \dp1=\dimen2 \box1\relax}

\def\wbb#1{\kern#1em}
\def\vci{\vrule  width.02em height1.47ex depth-.0ex}		
\def\11{{\rm\wbb{.2}\vci\wbb{-.37}1}}

\def\underset#1#2{\mathrel{\mathop{\kern0pt #2}\limits_{#1}}}

\def\overset#1#2{\mathrel{\mathop{\kern0pt #2}\limits^{#1}}}

\newtheorem{thm}{Theorem}[section]
\newtheorem{lem}[thm]{Lemma}
\newtheorem{prop}[thm]{Proposition}
\newtheorem{cor}[thm]{Corollary}
\newtheorem{defin}[thm]{Definition}
\newtheorem{rem}[thm]{Remark}

\begin{document}

\title{Sobolev solutions of parabolic equation in a complete riemannian manifold.}

\author{Eric Amar}

\date{}
\maketitle
 \renewcommand{\abstractname}{Abstract}

\begin{abstract}
We study Sobolev estimates for the solutions of parabolic equations
 acting on a vector bundle, in a complete riemannian manifold
 $M.$ The idea is to introduce geometric weights on $M.$ We get
 global Sobolev estimates with these weights. As applications,
 we find and improve "classical results", i.e. results without
 weights.  As an example we get Sobolev estimates for the solutions
 of the heat equation on $p$-forms when the manifold has "weak
 bounded geometry " of order $1$\!\!\!\! .\par 

\end{abstract}
\ \par 
\ \par 
\ \par 

\tableofcontents
\ \par 
\ \par 

\section{Introduction.}
\quad The study of $L^{r}$ estimates for the solutions of parabolic
 equations in a complete riemannian manifold started long time
 ago. For the case of the heat equation, a basic work was done
 by R.S. Strichartz~\cite{Strichartz83}. In particular he proved
 that the heat kernel is a contraction on the space of functions
 in $L^{r}(M)$ for $1\leq r\leq \infty .$ \ \par 
\quad Let $(M,g)$ be a complete riemannian manifold and let $G:=(H,\pi
 ,M)$ be a complex ${\mathcal{C}}^{m}$ vector bundle over $M$
 of rank $N$ with fiber $H.$\ \par 
Let $A$ be an elliptic operator of order $m$ acting on sections
 of $G$ to themselves.\ \par 
\quad Our aim here is to get Sobolev estimates on the solutions of
 the parabolic equation $Du:=\partial _{t}u-Au=\omega ,$ where
 $u,\omega $ are sections of $G$ over $M.$\ \par 
\quad Opposite to the usual way to do, see for instance the book by
 Grigor'yan~\cite{Grigor09} and the references therein or the
 paper by~\cite{Nistor06}, we do not use estimates on the kernel
 associated to the semi group of the differential operator on
 the manifold.\ \par 
\quad We shall follow another natural path to proceed: first we use
  known result in ${\mathbb{R}}^{n}$ to get precise local estimates
 on $M,$ then we globalise them.\ \par 
\quad The advantage of this way is that, for instance when dealing
 with the heat equation, we need no assumptions on the heat kernel.\ \par 
\ \par 
\quad To present the ideas in a simple way, we first restrict ourselves
 to the basic case of the heat equation $Du:=\partial _{t}u+\Delta
 u=\omega ,$ where $\Delta :=dd^{*}+d^{*}d$ is the Hodge laplacian
 and $u,\omega $ belong to the vector bundle of differential $p$-forms.\ \par 
\quad We introduce $(m,\epsilon )$-admissible balls $B_{m,\epsilon
 }(x)$ in $(M,g).$ These balls are the ones defined in the work
 of Hebey and Herzlich~\cite{HebeyHerzlich97} but without asking
 for the harmonicity of the local coordinates.\ \par 
\quad Then we use a Theorem by Haller-Dintelmann, Heck and Hieber~\cite[Corollary
 3.2, p. 5]{HHH06} done in ${\mathbb{R}}^{n},$ to get precise
 local results on these $(m,\epsilon )$-admissible balls.\ \par 
For $x$ in $M,$ the radius $R_{m,\epsilon }(x)$ of the admissible
 ball $B_{m,\epsilon }(x)$ tells us how far from the euclidean
 geometry of ${\mathbb{R}}^{n}$ the manifold $(M,g)$ is near
 the point $x,$ and so it not surprising that our geometric weights
 are functions of these radius.\ \par 
\quad Finally we use an adapted Vitali covering to globalise the local
 results we got. \ \par 
\ \par 
\quad Let $W^{k,r}_{p}(M,w)$ be the space of $p$-forms on $M$ belonging
 in the Sobolev space $W^{k,r}(M,w)$ with the weight $w.$ The
 same way $L^{r}_{p}(M,w)$ is the space of $p$-forms on $M$ belonging
 in the Lebesgue space $L^{r}(M,w)$ with the weight $w.$\ \par 
This gives us the following theorem, written here in the case
 of the heat equation:\ \par 

\begin{thm}
Let $M$ be a connected complete $n$-dimensional ${\mathcal{C}}^{2}$
 riemannian manifold without boundary. Let  $\displaystyle Du:=\partial
 _{t}u+\Delta u$ be the heat operator acting on the bundle $\displaystyle
 \Lambda ^{p}(M)$ of $p$-forms on $M.$ Let:\par 
\quad \quad \quad $\displaystyle R(x)=R_{2,\epsilon }(x),\ w_{1}(x):=R(x)^{r\delta
 },\ w_{2}(x):=R(x)^{r\gamma },\ w_{3}(x):=R(x)^{r\beta },$\par 
where $\beta ,\gamma ,\delta $ are explicit constants. Then,
 for any $\alpha >0,\ r\geq 2,$ we have:\par 
\quad \quad \quad $\displaystyle \forall \omega \in L^{r}(\lbrack 0,T+\alpha \rbrack
 ,L^{r}_{p}(M,w_{3}))\cap L^{r}(\lbrack 0,T+\alpha \rbrack ,L^{2}_{p}(M)),\
 \exists u\in L^{r}(\lbrack 0,T\rbrack ,W^{2,r}_{p}(M,w_{2}))::Du=\omega ,$\par 
with\par 
$\displaystyle {\left\Vert{\partial _{t}u}\right\Vert}_{L^{r}(\lbrack
 0,T\rbrack ,L^{r}_{p}(M,w_{1}))}+{\left\Vert{u}\right\Vert}_{L^{r}(\lbrack
 0,T\rbrack ,W^{2,r}_{p}(M,w_{2}))}\leq c_{1}{\left\Vert{\omega
 }\right\Vert}_{L^{r}(\lbrack 0,T+\alpha \rbrack ,L^{r}_{p}(M,w_{3}))}+c_{2}{\left\Vert{\omega
 }\right\Vert}_{L^{r}(\lbrack 0,T+\alpha \rbrack ,L_{p}^{2}(M))}.$\par 
In the case of functions instead of $p$-forms we have the same
 estimates but with $R(x)=R_{1,\epsilon }(x)$ and the weights:\par 
\quad \quad \quad $\displaystyle w_{1}(x):=R(x)^{r\delta '},\ w_{2}(x):=R(x)^{r\gamma
 '},\ w_{3}(x):=R(x)^{r\beta '}.$
\end{thm}
\quad Because our admissible radius $R_{m,\epsilon }(x)$ is smaller
 than one, to forget the weights, i.e. to get "classical estimates",
 it suffices to have $\forall x\in M,\ R_{m,\epsilon }(x)\geq
 \delta >0.$ In order to get this, we shall use a nice theorem
 by Hebey and Herzlich~\cite[Corollary, p. 7] {HebeyHerzlich97}
 which warranty us that the radius of our admissible balls is
 uniformly bounded below.\ \par 
\ \par 
\quad We introduce a weakened notion of bounded geometry: in the classical
 definition we replace the curvature tensor by the Ricci one:\ \par 

\begin{defin}
A riemannian manifold $M$ has $k$-order {\bf weak bounded geometry} if:\par 
\quad $\bullet $  the injectivity radius $r_{inj}(x)$ at $x\in M$ is
 bounded below by some constant $\delta >0$ for any $\displaystyle x\in M$\par 
\quad $\bullet $  for $0\leq j\leq k,$ the covariant derivatives $\nabla
 ^{j}R_{c}$ of the Ricci curvature tensor are bounded in $L^{\infty }(M)$ norm.
\end{defin}
\quad Using this notion we get the following theorem, written here
 in the case of the heat equation:\ \par 

\begin{thm}
Let $M$ be a connected complete $n$-dimensional ${\mathcal{C}}^{2}$
 riemannian manifold without boundary. Let  $\displaystyle Du:=\partial
 _{t}u+\Delta u$ be the heat operator acting on the bundle $\displaystyle
 \Lambda ^{p}(M)$ of $p$-forms on $M.$ Suppose moreover that
 $(M,g)$ has $1$ order weak bounded geometry. Then\par 
\quad \quad \quad $\displaystyle \forall \omega \in L^{r}(\lbrack 0,T+\alpha \rbrack
 ,L^{r}_{p}(M))\cap L^{r}(\lbrack 0,T+\alpha \rbrack ,L^{2}_{p}(M)),\
 \exists u\in L^{r}(\lbrack 0,T\rbrack ,W^{2,r}_{p}(M))::Du=\omega ,$\par 
with:\par 
\quad \quad \quad $\displaystyle {\left\Vert{\partial _{t}u}\right\Vert}_{L^{r}(\lbrack
 0,T\rbrack ,L_{p}^{r}(M))}+{\left\Vert{u}\right\Vert}_{L^{r}(\lbrack
 0,T\rbrack ,W_{p}^{2,r}(M))}\leq c_{1}{\left\Vert{\omega }\right\Vert}_{L^{r}(\lbrack
 0,T+\alpha \rbrack ,L_{p}^{r}(M))}+c_{2}{\left\Vert{\omega }\right\Vert}_{L^{r}(\lbrack
 0,T+\alpha \rbrack ,L_{p}^{2}(M))}.$\par 
In the case of functions instead of $p$-forms we have the same
 estimates just supposing that $(M,g)$ has $0$ order weak bounded geometry.
\end{thm}
\quad Our method extends to the study of general parabolic equation
 of order $m$ acting on metric vector bundles. But even in the
 special case of the heat equation acting on $p$-forms, it gives
 some new insights. Let us compare with 3 papers using the heat
 kernel method. These papers give estimates on the solutions
 of the heat equation $Du=\omega $ for $u(t,x)$ with $t\in \lbrack
 0,T\rbrack $ fixed. On the other hand the solutions I get are
 in $L^{r}(\lbrack 0,T\rbrack ,\ W^{m,r}_{p}(M)).$\ \par 
\quad $\bullet $ Comparing with the result of Strichartz~\cite{Strichartz83}
 on functions, he has no condition at all to get $u(t,\cdot )\in
 L^{r}(M)$ for $\omega (t,\cdot )\in L^{r}(M)\cap L^{2}(M)$ for
 any $r\in \lbrack 1,\infty \rbrack .$\ \par 
\quad Here we get $u\in L^{r}(\lbrack 0,T\rbrack ,\ W^{m,r}(M))$ for
 $\omega \in L^{r}(\lbrack 0,T+\alpha \rbrack ,\ L^{r}(M))\cap
 L^{r}(\lbrack 0,T+\alpha \rbrack ,\ L^{2}(M)),$ at the price
 that $(M,g)$ has $0$ order weak bounded geometry. Moreover,
 by Theorem 8.7 in~\cite{SobAmar19}, the Sobolev embeddings are
 true in that case, hence $u\in W^{2,r}(M)\Rightarrow u\in L^{s}(M)$
 with $\frac{1}{s}=\frac{1}{r}-\frac{2}{n},$ and the result is
 improved also in the Lebesgue scale.\ \par 
\quad $\bullet $ The work by~\cite{Nistor06}, also using the kernel
 associated to the semi group of the differential operator acting
 on metric vector bundles, contains a wide range of precise results,
 among them Sobolev estimates for the solutions of the parabolic
 equation. This is done under geometrical hypotheses on the manifold,
 essentially: bounded geometry of any order.\ \par 
\quad Here we allow the order $m$ of the parabolic equation to be greater
 that $2$ and we need only that $(M,g)$ has $m-1$ order weak
 bounded geometry to get Sobolev estimates, but the price is
 that we have our solutions in $L^{r}(\lbrack 0,T\rbrack ,\ W^{m,r}_{p}(M)),$
 not in $W^{m,r}_{p}(M),$ for any $t\in \lbrack 0,T\rbrack .$\ \par 
\quad $\bullet $ Comparing to the result in~\cite[Theorem 1.2]{MagOuh17},
 the hypotheses they have are directly on the kernel and on the
 manifold: the heat kernel must satisfy a Gaussian upper bound,
 $M$ must satisfy a volume doubling condition, plus another condition
 on the negative part of the Ricci curvature. They get Lebesgue
 estimates on $p$-forms $u(t,\cdot )\in L_{p}^{r}(M)$ for $\omega
 (t,\cdot )\in L_{p}^{r}(M).$\ \par 
\quad Here again we need that $(M,g)$ has $1$ order weak bounded geometry
 to get Sobolev estimates, which are better than Lebesgue estimates,
 but in $L^{r}(\lbrack 0,T\rbrack ,\ W^{m,r}_{p}(M)).$\ \par 
\ \par 
\quad The proofs here are, of course, completely different than the
 proofs using kernels.\ \par 

\section{Notation, definitions and main results.}

\subsection{Admissible balls.}

\begin{defin}
Let~\label{mLIR25} $(M,g)$ be a riemannian manifold and $\displaystyle
 x\in M.$ We shall say that the geodesic ball $\displaystyle
 B(x,R)$ is $(m,\epsilon )$-{\bf admissible} if there is a chart
 $\displaystyle (B(x,R),\varphi )$ such that, with $\epsilon \in (0,1)$:\par 
\quad 1) $\displaystyle (1-\epsilon )\delta _{ij}\leq g_{ij}\leq (1+\epsilon
 )\delta _{ij}$ in $\displaystyle B(x,R)$ as bilinear forms,\par 
\quad 2) $\displaystyle \ \sum_{1\leq \left\vert{\beta }\right\vert
 \leq m}{R^{\left\vert{\beta }\right\vert }\sup \ _{i,j=1,...,n,\
 y\in B_{x}(R)}\left\vert{\partial ^{\beta }g_{ij}(y)}\right\vert
 }\leq \epsilon .$\par 
We shall denote ${\mathcal{A}}_{m}(\epsilon )$ the set of $(m,\epsilon
 )$-admissible balls.
\end{defin}

\begin{defin}
~\label{CL26}Let $\displaystyle x\in M,$ we set $R'(x)=\sup \
 \lbrace R>0::B(x,R)\in {\mathcal{A}}_{m}(\epsilon )\rbrace .$
 We shall say that $\displaystyle R_{\epsilon }(x):=\min \ (1,R'(x)/2)$
 is the $(m,\epsilon )$-{\bf admissible radius} at $\displaystyle x.$
\end{defin}

\begin{rem}
Let $\displaystyle x,y\in M.$ Suppose that $R'(x)>d_{g}(x,y),$
 where $d_{g}(x,y)$ is the riemannian distance between $x$ and
 $y.$ Consider the ball $B(y,\rho )$ of center $y$ and radius
 $\rho :=R'(x)-d_{g}(x,y).$ This ball is contained in $B(x,R'(x))$
 hence, by definition of $R'(x),$ we have that all the points
 in $B(y,\rho )$ verify the conditions 1) and 2) so, by definition
 of $R'(y),$ we have that $R'(y)\geq R'(x)-d_{g}(x,y).$ If $\displaystyle
 R'(x)\leq d_{g}(x,y)$ this is also true because $\displaystyle
 R'(y)>0.$ Exchanging $x$ and $y$ we get that $\left\vert{R'(y)-R'(x)}\right\vert
 \leq d_{g}(x,y).$\par 
\quad Hence $R'(x)$ is $1$-lipschitzian so it is continuous. So the
 $\epsilon $-admissible radius $R_{\epsilon }(x)$ is continuous.
\end{rem}

\begin{rem}
~\label{SC16} Because an admissible ball $B(x,R_{\epsilon }(x))$
 is geodesic, we get that the injectivity radius $r_{inj}(x)$
 always verifies $\displaystyle r_{inj}(x)\geq R_{\epsilon }(x).$
\end{rem}

\begin{lem}
~\label{m7}(Slow variation of the admissible radius) Let $(M,g)$
 be a riemannian manifold. With $R(x)=R_{\epsilon }(x)=$ the
 $\epsilon $-admissible radius at $x\in M,\ \forall y\in B(x,R(x))$
 we have $R(x)/2\leq R(y)\leq 2R(x).$
\end{lem}
\quad Proof.\ \par 
Let $x,y\in M$ and $d(x,y)$ the riemannian distance on $(M,g).$
 Let $y\in B(x,R(x))$ then $d(x,y)\leq R(x)$ and suppose first
 that $R(x)\geq R(y).$ \ \par 
Then, because $R(x)=R'(x)/2,$ we get $y\in B(x,R'(x)/2)$ hence
 we have $B(y,R'(x)/2)\subset B(x,R'(x)).$ But by the definition
 of $R'(x),$ the ball $\displaystyle B(x,R'(x))$ is admissible
 and this implies that the ball $\displaystyle B(y,R'(x)/2)$
 is also admissible for exactly the same constants and the same
 chart; this implies that $R'(y)\geq R'(x)/2$ hence $R(y)\geq
 R(x)/2,$ so $R(x)\geq R(y)\geq R(x)/2.$\ \par 
\quad If $R(x)\leq R(y)$ then\ \par 
\quad \quad \quad $\displaystyle d(x,y)\leq R(x)\Rightarrow d(x,y)\leq R(y)\Rightarrow
 x\in B(y,R'(y)/2)\Rightarrow B(x,R'(y)/2)\subset B(y,R'(y)).$\ \par 
Hence the same way as above we get $R(y)\geq R(x)\geq R(y)/2\Rightarrow
 R(y)\leq 2R(x).$ So in any case we proved that\ \par 
\quad \quad \quad $\displaystyle \forall y\in B(x,R(x)),\ R(x)/2\leq R(y)\leq 2R(x).$
  $\blacksquare $\ \par 

\subsection{Vector bundle.}
\quad Let $(M,g)$ be a complete riemannian manifold and let $G:=(H,\pi
 ,M)$ be a complex ${\mathcal{C}}^{m}$ vector bundle over $M$
 of rank $N$ with fiber $H.$ Suppose moreover that $G$ has a
 smooth scalar product $(\ ,\ )$ and  a \emph{metric} connection
 $\nabla ^{G}:{\mathcal{C}}^{\infty }(M,G)\rightarrow {\mathcal{C}}^{\infty
 }(M,G\otimes T^{*}M),$ i.e. verifying $\displaystyle d(u,v)=(\nabla
 ^{G}u,v)+(u,\nabla ^{G}v),$ where $d$ is the exterior derivative
 on $M$ acting on the scalar product $(u,v).$ See ~\cite[Section
 13]{TaylorGD}.\ \par 

\begin{lem}
~\label{SB30} The $\epsilon $-admissible balls $B(x,R_{\epsilon
 }(x))$ trivialise the bundle $G.$
\end{lem}
\quad Proof.\ \par 
Because if $B(x,R)$ is a $\epsilon $-admissible ball, we have
 by Remark~\ref{SC16} that $R\leq r_{inj}(x).$ Then, one can
 choose a local frame field for $G$ on $\displaystyle B(x,R)$
 by radial parallel translation, as done in~\cite[Section 13,
 p.86-87]{TaylorGD}, see also~\cite[p. 4, eq. (1.3)]{Nistor06}.
 This means that the $\epsilon $-admissible ball also trivialises
 the bundle $G.$ $\blacksquare $\ \par 
\ \par 
\quad If $\partial _{j}:=\partial /\partial x_{j}$ in a coordinate
 system on, say $B(x_{0},R),$ and with a local frame $\lbrace
 e_{\alpha }\rbrace _{\alpha =1,...,N},$ we have, for a smooth
 sections of $G,\ u=u^{\alpha }e_{\alpha }$ with the Einstein
 summation convention. We set:\ \par 
\quad \quad \quad $\displaystyle \nabla _{\partial _{j}}u=(\partial _{j}u^{\alpha
 }\ +u^{\beta }\Gamma ^{G,\alpha }_{\beta j})e_{\alpha },$\ \par 
the Christoffel coefficients $\displaystyle \Gamma ^{G,\alpha
 }_{\beta j}$ being defined by $\displaystyle \nabla _{\partial
 _{j}}e_{\beta }=\Gamma ^{G,\alpha }_{\beta j}e_{\alpha }.$\ \par 
\quad We shall make the following hypothesis on the connection on $G,$
 for $B(x_{0},R)\in {\mathcal{A}}(m,\epsilon )$:\ \par 
(CMT)        $\displaystyle \forall x\in B(x_{0},R),\ \forall
 k\leq \ m,\ \ \left\vert{\partial ^{k-1}\Gamma ^{G,\alpha }_{\beta
 j}(x)}\right\vert \leq C(n,G,\epsilon )\sum_{\left\vert{\beta
 }\right\vert \leq k}{\sup \ _{i,j=1,...,n,}\left\vert{\partial
 ^{\beta }g_{ij}(x)}\right\vert },$\ \par 
the constant $C$ depending only on $n,\epsilon $ and $G$ but
 not on $\displaystyle B(x_{0},R)\in {\mathcal{A}}(m,\epsilon )$\ \par 
\quad This hypothesis is natural:\ \par 

\begin{lem}
~\label{pBB26}The hypothesis (CMT) is true for the Levi-Civita
 connection on $M.$
\end{lem}
\quad Proof.\ \par 
Let $\displaystyle \Gamma ^{k}_{lj}$ be the Christoffel coefficients
 of the Levi-Civita connection  on the tangent bundle $TM.$ We have\ \par 
\quad \quad \quad \begin{equation} \Gamma ^{i}_{kj}=\frac{1}{2}g^{il}(\frac{\partial
 g_{kl}}{\partial x^{j}}+\frac{\partial g_{lj}}{\partial x^{k}}-\frac{\partial
 g_{jk}}{\partial x^{l}}).\label{SB34}\end{equation}\ \par 
Now on $\displaystyle B(x_{0},R)\in {\mathcal{A}}_{m}(\epsilon
 ),$ we have $\displaystyle (1-\epsilon )\delta _{ij}\leq g_{ij}\leq
 (1+\epsilon )\delta _{ij}$ as bilinear forms. Hence\ \par 
\quad \quad \quad $\displaystyle \forall x\in B(x_{0},\ R),\ \ \left\vert{\Gamma
 ^{i}_{kj}(x)}\right\vert \leq \frac{3}{2}(1-\epsilon )^{-1}\sum_{\left\vert{\beta
 }\right\vert =1}{\sup \ _{i,j=1,...,n,}\left\vert{\partial ^{\beta
 }g_{ij}(x)}\right\vert }$\ \par 
in a coordinates chart on $\displaystyle B(x_{0},R).$\ \par 
\quad We have the same with~(\ref{SB34}) for the derivatives of $\displaystyle
 \Gamma ^{i}_{kj}.$ $\blacksquare $\ \par 

\begin{rem}
~\label{pBB24} If the hypothesis (CMT) is true for two vector
 bundles on $M,$ then a easy computation gives that it is true
 for the tensor product of the two bundles over $M.$ In particular
 (CMT) is true for tensor bundles over $M.$ It is also true for
 the sub-bundle of $p$-forms on $M.$
\end{rem}

\subsection{Sobolev spaces for sections of $G$ with weight.~\label{p25}}
\quad We have seen that $\displaystyle \nabla ^{G}:\ {\mathcal{C}}^{\infty
 }(M,G)\rightarrow {\mathcal{C}}^{\infty }(M,G\otimes T^{*}M).$
 On the tensor product of two Hilbert spaces we put the canonical
 scalar product $(u\otimes \omega ,v\otimes \mu ):=(u,v)(\omega
 ,\mu ),$ with $u\otimes \omega \in G\otimes T^{*}M,$ and completed
 by linearity to all elements of the tensor product. On $T^{*}M$
 we have the Levi-Civita connection  $\nabla ^{M},$ which is
 of course a metric one, and on $G$ we have the metric connection
  $\nabla ^{G}$ so we define a connection  on the tensor product
 $\displaystyle G\otimes T^{*}M$:\ \par 
\quad \quad \quad $\displaystyle \nabla ^{G\otimes T^{*}M}(u\otimes \omega )=(\nabla
 ^{G}u)\otimes \omega +u\otimes (\nabla ^{T^{*}M}\omega )$\ \par 
by asking that this connection  be a derivation. We get easily that\ \par 
\quad \quad \quad $\displaystyle \nabla ^{G\otimes T^{*}M}\ :\ {\mathcal{C}}^{\infty
 }(M,G\otimes T^{*}M)\rightarrow {\mathcal{C}}^{\infty }(M,G\otimes
 (T^{*}M)^{\otimes 2})$\ \par 
is still a \emph{metric} connection , i.e.\ \par 
\quad \quad \quad $\displaystyle d(u\otimes \omega ,v\otimes \mu )=(\nabla ^{G\otimes
 T^{*}M}(u\otimes \omega ),\ v\ \otimes \mu )+(u\otimes \omega
 ,\ \nabla ^{G\otimes T^{*}M}(v\otimes \mu )).$\ \par 
\quad We define by iteration $\nabla ^{j}u:=\nabla (\nabla ^{j-1}u)$
 on the section $u$ of $G$ and the associated pointwise scalar
 product $\displaystyle (\nabla ^{j}u(x),\nabla ^{j}v(x))$ which
 is defined on $\displaystyle G\otimes (T^{*}M)^{\otimes j},$
 with again the metric connection\ \par 
\quad \quad \quad $\displaystyle d(\nabla ^{j}u,\nabla ^{j}v)(x)=(\nabla ^{j+1}u,\nabla
 ^{j}v)(x)+(\nabla ^{j}u,\nabla ^{j+1}v)(x).$\ \par 
\quad Let $w$ be a weight on $M,$ i.e. a positive measurable function
 on $M.$ If $\displaystyle k\in {\mathbb{N}}$ and $r\geq 1$ are
 given, we denote by ${\mathcal{C}}^{k,r}_{G}(M,w)$ the space
 of smooth sections of $G$ $\omega \in {\mathcal{C}}^{\infty
 }(M)$ such that $\ \left\vert{\nabla ^{j}\omega }\right\vert
 \in L^{r}(M,w)$ for $j=0,...,k$ with the pointwise modulus associated
 to the pointwise scalar product. Hence\ \par 
\quad \quad \quad $\displaystyle {\mathcal{C}}^{k,r}_{G}(M,w):=\lbrace \omega \in
 {\mathcal{C}}_{G}^{\infty }(M),\ \forall j=0,...,k,\ \int_{M}{\left\vert{\nabla
 ^{j}\omega }\right\vert ^{r}(x)w(x)dv(x)}<\infty \rbrace ,$\ \par 
with $dv$ the volume measure on $(M,g).$\ \par 
\quad Now we have, see M. Cantor ~\cite[Definition 1 \&  2, p. 240]{Cantor74}
 for the case without weight:\ \par 

\begin{defin}
The Sobolev space $\displaystyle W_{G}^{k,r}(M,w)$ is the completion
 of ${\mathcal{C}}^{k,r}_{G}(M,w)$ with respect to the norm:\par 
\quad \quad \quad $\displaystyle {\left\Vert{\omega }\right\Vert}_{W_{G}^{k,r}(M,w)}=\sum_{j=0}^{k}{{\left({\int_{M}{\left\vert{\nabla
 ^{j}\omega (x)}\right\vert ^{r}w(x)dv(x)}}\right)}^{1/r}}.$
\end{defin}
The usual case is when $w\equiv 1.$ Then we write simply $\displaystyle
 W_{G}^{k,r}(M).$\ \par 
\ \par 
\quad A vector bundle $G$ verifying the following two hypotheses will
 be called \textbf{adapted}:\ \par 
\quad $\bullet $ the vector bundle $G$ is equipped with  a \emph{metric
 connection};\ \par 
\quad $\bullet $ the Christoffel symbols $\Gamma ^{G,\alpha }_{\beta
 j}$ of the connection are controlled by the metric tensor (CMT) $g$:\ \par 
(CMT)        $\forall x\in B(x_{0},R),\ \forall k\leq m,\ \ \left\vert{\partial
 ^{k-1}\Gamma ^{G,\alpha }_{\beta j}(x)}\right\vert \leq C(n,G,\epsilon
 )\sum_{\left\vert{\beta }\right\vert \leq k}{\sup \ _{i,j=1,...,n,}\left\vert{\partial
 ^{\beta }g_{ij}(x)}\right\vert },$\ \par 
the constant $C$ depending only on $n,\epsilon $ and $G$ but
 not on the admissible ball $B(x_{0},R)\in {\mathcal{A}}_{m}(\epsilon ).$\ \par 

\subsection{Parabolic operator.}
\quad We suppose that $Du:=\partial _{t}u-Au$ is parabolic in ${\mathbb{R}}^{n}$
 in the sense of~\cite{HHH06}:\ \par 
\quad $\bullet $ $A$ is a system of differential operators of the form
 $A=\sum_{\left\vert{\alpha }\right\vert \leq m}{a_{\alpha }\partial
 ^{\alpha }},$ where $\partial =-i(\partial _{1},...,\partial
 _{n})$ and $a_{\alpha }\in L^{\infty }({\mathbb{R}}^{n},{\mathbb{C}}^{N{\times}N}).$\
 \par 
\quad $\bullet $ $A$ is $(C,\theta )$-elliptic; this means that there
 exist constants $\theta \in \lbrack 0,\pi )$ and $C>0,$ such
 that the principal part $A_{\# }(x,\xi ):=\sum_{\left\vert{\alpha
 }\right\vert =m}{a_{\alpha }\xi ^{\alpha }}$ of the symbol of
 $A$ satisfies the following conditions:\ \par 
\quad \quad \quad $\sigma (A_{\# }(x,\xi ))\subset \bar S_{\theta }$ and ${\left\Vert{A_{\#
 }(x,\xi )^{-1}}\right\Vert}\leq M$ for all $\xi \in {\mathbb{R}}^{n},\
 \left\vert{\xi }\right\vert =1,$\ \par 
for almost all $x\in {\mathbb{R}}^{n}.$ Here $\displaystyle S_{\theta
 }$ denotes the sector in the complex plane defined by $\displaystyle
 S_{\theta }:=\lbrace \lambda \in {\mathbb{C}}\backslash \lbrace
 0\rbrace ::\left\vert{\mathrm{a}\mathrm{r}\mathrm{g}\lambda
 }\right\vert <\theta \rbrace $ and the spectrum of an $N{\times}N$-matrix
 ${\mathcal{M}}$ is denoted by $\sigma ({\mathcal{M}}).$\ \par 
\quad $\bullet $ Because we work only with the usual Lebesgue spaces,
 we take for the domain of $A,\ {\mathcal{D}}(A):=W^{m,r}({\mathbb{R}}^{n})^{N}.$\
 \par 
We shall use the following~\cite[Corollary 3.2, p. 5]{HHH06}:\ \par 

\begin{thm}
~\label{pLIR0}Let $n\geq 2,\ 1<r,s<\infty ,\ \theta \in (0,\pi
 )$ and $C>0.$ Assume that $A:=\sum_{\left\vert{\alpha }\right\vert
 \leq m}{a_{\alpha }(x)\partial ^{\alpha }}$ is a $(C,\theta
 )$-elliptic operator in $L^{r}_{w}({\mathbb{R}}^{n})^{N}$ with
 coefficients $a_{\alpha }$ satisfying:\par 
\quad a) $a_{\alpha }\in L^{\infty }({\mathbb{R}}^{n};{\mathbb{C}}^{N{\times}N})\cap
 VMO({\mathbb{R}}^{n};{\mathbb{C}}^{N{\times}N})$ for $\left\vert{\alpha
 }\right\vert =m,$\par 
\quad b) $\displaystyle a_{\alpha }\in L^{\infty }({\mathbb{R}}^{n};{\mathbb{C}}^{N{\times}N})$
 for $\left\vert{\alpha }\right\vert <m.$\par 
Suppose that $\displaystyle Du:=\partial _{t}u-Au=\omega ,\ u(x,0)\equiv
 0,$ and assume now that $\theta <\frac{\pi }{2},$ then there
 exist constants $M,\mu \geq 0$ such that, with $J:=\lbrack 0,\infty
 \lbrack ,$\par 
\quad \quad \quad $\displaystyle {\left\Vert{\partial _{t}u}\right\Vert}_{L^{s}(J,L^{r}({\mathbb{R}}^{n})^{N})}+{\left\Vert{(\mu
 +A)u}\right\Vert}_{L^{s}(J,L^{r}({\mathbb{R}}^{n})^{N})}\leq
 M{\left\Vert{\omega }\right\Vert}_{L^{s}(J,L^{r}({\mathbb{R}}^{n})^{N})}.$\par 
Moreover the solution $u$ is unique verifying this estimate.
\end{thm}

\subsection{Global assumptions.}
\quad We shall made the following global assumption on the operator
 $A$ in the riemannian manifold $M$ in all the sequel of this work.\ \par 

\begin{defin}
We say that the operator $A$ is $(C,\theta )$-elliptic of order
 $m$ acting on sections of  $G$ in the riemannian manifold $(M,g),$
 if for any chart $(U,\varphi )$ on $(M,g)$ which trivializes
 $G,$ i.e. $G_{\varphi },$ the image of $G,$ is the trivial bundle
  $\varphi (U){\times}{\mathbb{R}}^{N}$ in $\displaystyle \varphi
 (U),$ we have, with $A_{\varphi }$ the image of the operator $A$:\par 
\quad $\bullet $ $A_{\varphi }$ is a system of differential operators
 of the form $A_{\varphi }=\sum_{\left\vert{\alpha }\right\vert
 \leq m}{a_{\alpha }\partial ^{\alpha }},$ where $\partial =-i(\partial
 _{1},...,\partial _{n})$ and $a_{\alpha }\in L^{\infty }(\varphi
 (U),{\mathbb{C}}^{N{\times}N}),$ with:\par 
\quad a) $a_{\alpha }\in L^{\infty }(\varphi (U);{\mathbb{C}}^{N{\times}N})\cap
 VMO(\varphi (U);{\mathbb{C}}^{N{\times}N})$ for $\left\vert{\alpha
 }\right\vert =m,$\par 
\quad b) $\displaystyle a_{\alpha }\in L^{\infty }(\varphi (U);{\mathbb{C}}^{N{\times}N})$
 for $\left\vert{\alpha }\right\vert <m.$\par 
\par 
\quad $\bullet $ $A_{\varphi }$ is $(C,\theta )$-elliptic; this means
 that there exist constants $\theta \in \lbrack 0,\pi )$ and
 $C>0,$ such that the principal part $A_{\# }(x,\xi ):=\sum_{\left\vert{\alpha
 }\right\vert =m}{a_{\alpha }\xi ^{\alpha }}$ of the symbol of
 $A$ satisfies the following conditions:\S \par 
\quad \quad \quad $\sigma (A_{\# }(x,\xi ))\subset \bar S_{\theta }$ and ${\left\Vert{A_{\#
 }(x,\xi )^{-1}}\right\Vert}\leq M$ for all $\xi \in {\mathbb{R}}^{n},\
 \left\vert{\xi }\right\vert =1,$\par 
for almost all $x\in \varphi (U)^{n}.$ \par 
And all the bounds being independent of the chart $(U,\varphi ).$
\end{defin}
\ \par 
\quad We shall also need the following "\textbf{threshold hypothesis}".\ \par 
(THL2)  For any $\omega \in L^{s}(\lbrack 0,T\rbrack ,L^{2}_{G}(M))$
 there is a $\displaystyle u\in L^{s}(\lbrack 0,T\rbrack ,L^{2}_{G}(M))$
 such that $Du=\omega $ with the estimate:\ \par 
\quad \quad \quad $\displaystyle {\left\Vert{u}\right\Vert}_{L^{s}(\lbrack 0,T\rbrack
 ,L^{2}_{G}(M))}\lesssim {\left\Vert{\omega }\right\Vert}_{L^{s}(\lbrack
 0,T\rbrack ,L^{2}_{G}(M))}.$\ \par 
This hypothesis is natural in the sense that it is true for the
 heat equation.\ \par 

\subsection{Main results.}
\quad We shall use the following notation.\ \par 

\begin{defin}
~\label{pBB29} For $r\geq 2,\ m\in {\mathbb{N}},\ m\geq 1,$ let
 $k:={\left\lceil{\frac{n(r-2)}{2mr}}\right\rceil}$ and define:\par 
\quad if $k=0,\ \beta (r,m):=m+\frac{n}{2}-\frac{n}{r};$\par 
\quad if $k\geq 1,$\par 
\quad \quad \quad $\displaystyle \beta =\beta (r,m):=\min (m+\frac{n}{2}-\frac{n}{r},\
 5m);\ \gamma =\gamma (r,m)=(4k+2)m;\ \delta =\delta (r,m)=(4k+1)m.$\par 
Define also:\par 
\quad if $k=0,\ \beta '=\beta '(r,m):=m+\frac{n}{2}-\frac{n}{r};$\par 
\quad if $k\geq 1,$\par 
\quad \quad \quad $\displaystyle \beta '=\beta '(r,m):=\min (m+\frac{n}{2}-\frac{n}{r},\
 4m);\ \gamma '=\gamma '(r,m)=(4m-1)k\ +2m;\ \delta '=\delta '(r,m)=(4m-1)k+m.$
\end{defin}
\quad We are in position to state the first main result of this work.\ \par 

\begin{thm}
Let $M$ be a connected complete $n$-dimensional ${\mathcal{C}}^{m}$
 riemannian manifold without boundary. Let $G:=(H,\pi ,M)$ be
 a complex ${\mathcal{C}}^{m}$ adapted vector bundle over $M.$
 Suppose $\displaystyle Du:=\partial _{t}u-Au,$ where $A$ is
 $(C,\theta )$-elliptic of order $m$ acting on sections of  $G$
 with $\theta <\pi /2$ in $(M,g).$ Moreover suppose we have (THL2).
 Let $r\geq 2$ and:\par 
\quad \quad \quad $\displaystyle R(x)=R_{m,\epsilon }(x),\ w_{1}(x):=R(x)^{r\delta
 },\ w_{2}(x):=R(x)^{r\gamma },\ w_{3}(x):=R(x)^{r\beta },$\par 
with $\beta ,\gamma ,\delta $ as in Definition~\ref{pBB29}. Then,
 for any $\alpha >0,\ r\geq 2,$ we have:\par 
\quad \quad \quad $\displaystyle \forall \omega \in L^{r}(\lbrack 0,T+\alpha \rbrack
 ,L^{r}_{G}(M,w_{3}))\cap L^{r}(\lbrack 0,T+\alpha \rbrack ,L^{2}_{G}(M)),\
 \exists u\in L^{r}(\lbrack 0,T\rbrack ,W^{m,r}_{G}(M))::Du=\omega ,$\par 
with\par 
\quad \quad \quad $\displaystyle {\left\Vert{\partial _{t}u}\right\Vert}_{L^{r}(\lbrack
 0,T\rbrack ,L^{r}_{G}(M,w_{1}))}+{\left\Vert{u}\right\Vert}_{L^{r}(\lbrack
 0,T\rbrack ,W^{m,r}_{G}(M,w_{2}))}\leq $\par 
\quad \quad \quad \quad \quad \quad \quad $\displaystyle \leq c_{1}{\left\Vert{\omega }\right\Vert}_{L^{r}(\lbrack
 0,T+\alpha \rbrack ,L^{r}_{G}(M,w_{3}))}+c_{2}{\left\Vert{\omega
 }\right\Vert}_{L^{r}(\lbrack 0,T+\alpha \rbrack ,L_{G}^{2}(M))},$\par 
In the case of functions instead of sections of $G$ we have the
 same estimates but with $R(x)=R_{m-1,\epsilon }(x)$ and the weights:\par 
\quad \quad \quad $\displaystyle w_{1}(x):=R(x)^{r\delta '},\ w_{2}(x):=R(x)^{r\gamma
 '},\ w_{3}(x):=R(x)^{r\beta '}.$
\end{thm}
\quad To find and improve "classical results", i.e. results without
 weights, we use a Theorem by Hebey and Herzlich~\cite[Corollary,
 p. 7] {HebeyHerzlich97} which warranty us that the radius of
 our "admissible balls" is uniformly bounded below.\ \par 
\quad This gives the second main result of this work.\ \par 

\begin{thm}
Suppose that $A$ is a $(C,\theta )$-elliptic operator of order
 $m$ acting on sections of the adapted vector bundle $\displaystyle
 G:=(H,\pi ,M)$ in the complete riemannian manifold $(M,g),$
 with $\theta <\pi /2.$ Consider the parabolic equation $\displaystyle
 Du=\partial _{t}u-Au$ also acting on sections of $G.$ Suppose
 moreover that $(M,g)$ has $(m-1)$ order weak bounded geometry
 and (THL2) is true. Let $r\geq 2$ then\par 
\quad \quad \quad $\displaystyle \forall \omega \in L^{r}(\lbrack 0,T+\alpha \rbrack
 ,L^{r}_{G}(M))\cap L^{r}(\lbrack 0,T+\alpha \rbrack ,L^{2}_{G}(M)),\
 \exists u\in L^{r}(\lbrack 0,T\rbrack ,W^{m,r}_{G}(M))::Du=\omega ,$\par 
with:\par 
\quad $\displaystyle {\left\Vert{\partial _{t}u}\right\Vert}_{L^{r}(\lbrack
 0,T\rbrack ,L_{G}^{r}(M))}+{\left\Vert{u}\right\Vert}_{L^{r}(\lbrack
 0,T\rbrack ,W_{G}^{m,r}(M))}\leq c_{1}{\left\Vert{\omega }\right\Vert}_{L^{r}(\lbrack
 0,T+\alpha \rbrack ,L_{G}^{r}(M))}+c_{2}{\left\Vert{\omega }\right\Vert}_{L^{r}(\lbrack
 0,T+\alpha \rbrack ,L_{G}^{2}(M))}.$\par 
In the case of functions instead of sections of $G$ we have the
 same estimates just supposing that $(M,g)$ has $(m-2)$ order
 weak bounded geometry.
\end{thm}

\section{Local results.}

\subsection{Local results in ${\mathbb{R}}^{n}.$}
\quad The following result follows the lines of ~\cite[Theorem 3.5]{ellipticEq18}:\
 \par 

\begin{thm}
~\label{m0}Let $A$ be an operator of order $m$ on $G$ in the
 complete riemannian manifold $M.$ Suppose that $A$ is elliptic
 and with ${\mathcal{C}}^{1}(M)$ smooth coefficients. Then, for
 any $x\in M$ and any ball $B:=B(x,R)$ such that $B(x,R)$ is
 a basis of a chart of $M$ around $x$ and trivialises the bundle
 $G,$ with the ball $\displaystyle B^{1}:=B(x,R/2),$ we have:\par 
\quad \quad \quad $\displaystyle {\left\Vert{u}\right\Vert}_{W^{m,r}_{G}(B^{1})}\leq
 c_{1}{\left\Vert{Au}\right\Vert}_{L^{r}_{G}(B)}+c_{2}R^{-m}{\left\Vert{u}\right\Vert}_{L^{r}_{G}(B)}.$\par
 
Moreover the constants are independent of the radius $R$ of the ball $B.$
\end{thm}
\quad We shall use Theorem~\ref{m0} in the proof of the following precise
 interior regularity theorem in the case of ${\mathbb{R}}^{n}.$
  The point here is that we need to have a clear dependence in
 the radius $R.$\ \par 

\begin{thm}
~\label{m1}Suppose that $A$ is a system of differential operators
 $(C,\theta )$-elliptic with $\theta <\pi /2,$ operating in ${\mathbb{R}}^{n},$
 and suppose $u$ is any solution of  the parabolic equation $\displaystyle
 Du=\partial _{t}u-Au=\omega $ in a ball $B(0,R)$ with $\displaystyle
 \omega \in L^{s}(\lbrack 0,T+\alpha \rbrack ,L^{r}(B)^{N})$
 and $\displaystyle u\in L^{s}(\lbrack 0,T+\alpha \rbrack ,L^{r}(B)^{N}).$\par 
Consider the ball $B^{1}:=B(0,R/2).$ We have, with $\alpha >0,\
 T>0$ and $r,s$ in $(1,\infty )$:\par 
\quad \quad \quad $\displaystyle {\left\Vert{\partial _{t}(u)}\right\Vert}_{L^{s}(\lbrack
 0,T\rbrack ,L^{r}(B^{1})^{N})}+{\left\Vert{u}\right\Vert}_{L^{s}(\lbrack
 0,T\rbrack ,W^{m,r}(B^{1})^{N})}\leq $\par 
\quad \quad \quad \quad \quad \quad \quad $\displaystyle \leq c_{1}{\left\Vert{D(u)}\right\Vert}_{L^{s}(\lbrack
 0,T+\alpha \rbrack ,L^{r}(B)^{N})}+c_{2}R^{-m}{\left\Vert{u}\right\Vert}_{L^{s}(\lbrack
 0,T+\alpha \rbrack ,L^{r}(B)^{N})}.$\par 
the constants $c_{j}$ being independent of $R.$
\end{thm}
\quad Proof.\ \par 
Let $\chi \in {\mathcal{D}}(B)$ such that $\chi (x)=1$ for  $x\in
 B^{1}.$ To ease the notation, let us set $\displaystyle L(s,r):=L^{s}(J,L^{r}({\mathbb{R}}^{n})^{N}).$\
 \par 
Because $A$ is  $(C,\theta )$-elliptic we can use the uniqueness
 in Theorem~\ref{pLIR0} to get that $v:=\chi u$ is the unique
 solution of $D(v)=D(\chi u)$ verifying, with $c_{1}$ independent of $B,$\ \par 
\quad \quad \quad $\displaystyle {\left\Vert{\partial _{t}(\chi u)}\right\Vert}_{L(s,r)}+{\left\Vert{(\mu
 +A)(\chi u)}\right\Vert}_{L(s,r)}\leq c_{1}{\left\Vert{D(\chi
 u)}\right\Vert}_{L(s,r)}.$\ \par 
Because\ \par 
\quad \quad \quad $\displaystyle {\left\Vert{(\mu +A)(\chi u)}\right\Vert}_{L(s,r)}\geq
 {\left\Vert{A(\chi u)}\right\Vert}_{L(s,r)}-\mu {\left\Vert{(\chi
 u)}\right\Vert}_{L(s,r)}$\ \par 
we have:\ \par 
\quad \quad \quad \begin{equation} {\left\Vert{\partial _{t}(\chi u)}\right\Vert}_{L(s,r)}+{\left\Vert{A(\chi
 u)}\right\Vert}_{L(s,r)}\leq c_{1}{\left\Vert{D(\chi u)}\right\Vert}_{L(s,r)}+\mu
 {\left\Vert{(\chi u)}\right\Vert}_{L(s,r)}.\label{pBB21}\end{equation}\ \par 
\quad We shall now use the estimates given by the ellipticity of $A.$
 For $t$ fixed, we have, by Theorem~\ref{m0}:\ \par 
\quad \quad \quad $\displaystyle {\left\Vert{\chi u}\right\Vert}_{W^{m,r}({\mathbb{R}}^{n})^{N}}\leq
 c_{2}{\left\Vert{A(\chi u)}\right\Vert}_{L^{r}({\mathbb{R}}^{n})^{N}}+c_{3}R^{-m}{\left\Vert{\chi
 u}\right\Vert}_{L^{r}({\mathbb{R}}^{n})^{N}}$\ \par 
where $c_{2},c_{3}$ are independent of $R.$\ \par 
So we get, integrating in $t$ and setting $\displaystyle W(s,r):=L^{s}(J,W^{m,r}({\mathbb{R}}^{n})^{N}),$\
 \par 
\quad \quad \quad $\displaystyle {\left\Vert{\chi u}\right\Vert}_{W(s,r)}\leq c_{2}{\left\Vert{A(\chi
 u)}\right\Vert}_{L(s,r)}+c_{3}R^{-m}{\left\Vert{\chi u}\right\Vert}_{L(s,r)}.$\
 \par 
Hence\ \par 
\quad \quad \quad $\displaystyle {\left\Vert{\partial _{t}(\chi u)}\right\Vert}_{L(s,r)}+{\left\Vert{\chi
 u}\right\Vert}_{W(s,r)}\leq {\left\Vert{\partial _{t}(\chi u)}\right\Vert}_{L(s,r)}+c_{2}{\left\Vert{A(\chi
 u)}\right\Vert}_{L(s,r)}+c_{3}R^{-m}{\left\Vert{\chi u}\right\Vert}_{L(s,r)}.$\
 \par 
Putting this in~(\ref{pBB21}) we get with $c_{4}:=\max (1,c_{2})$:\ \par 
\quad \quad \quad $\displaystyle {\left\Vert{\partial _{t}(\chi u)}\right\Vert}_{L(s,r)}+{\left\Vert{\chi
 u}\right\Vert}_{W(s,r)}\leq c_{4}c_{1}{\left\Vert{D(\chi u)}\right\Vert}_{L(s,r)}+c_{4}\mu
 {\left\Vert{(\chi u)}\right\Vert}_{L(s,r)}+c_{3}R^{-m}{\left\Vert{\chi
 u}\right\Vert}_{L(s,r)}.$\ \par 
So with new constants depending on $c_{1},c_{2},\ c_{3}$ and
 $\mu $ only and with $R\leq 1,$ we get\ \par 
\quad \quad \quad \begin{equation} {\left\Vert{\partial _{t}(\chi u)}\right\Vert}_{L(s,r)}+{\left\Vert{\chi
 u}\right\Vert}_{W(s,r)}\leq c'_{1}{\left\Vert{D(\chi u)}\right\Vert}_{L(s,r)}+c'_{2}R^{-m}{\left\Vert{\chi
 u}\right\Vert}_{L(s,r)}.\label{pBB23}\end{equation}\ \par 
Now we want to control $\displaystyle {\left\Vert{D(\chi u)}\right\Vert}_{L(s,r)}$
 by $\displaystyle {\left\Vert{D(u)}\right\Vert}_{L(s,r)}.$ We
 have, because $\chi $ does not depend on $t$:\ \par 
\quad \quad \quad \begin{equation} D(\chi u)=\chi \partial _{t}u-\chi Au+E=\chi
 Du+E\label{pBB22}\end{equation}\ \par 
with $E:=\chi Au-A(\chi u).$ The point is that $E$ contains only
 derivatives of the $j^{th}$ component of $u$ of order strictly
 less than in the $j^{th}$ component of $u$ in $Du.$\ \par 
So we have, fixing $t,$\ \par 
\quad \quad \quad $\displaystyle {\left\Vert{E}\right\Vert}_{L^{r}({\mathbb{R}}^{n})^{N}}\leq
 {\left\Vert{\partial \chi }\right\Vert}_{\infty }{\left\Vert{\chi
 u}\right\Vert}_{W^{m-1,r}({\mathbb{R}}^{n})^{N}}\leq R^{-1}{\left\Vert{\chi
 u}\right\Vert}_{W^{m-1,r}({\mathbb{R}}^{n})^{N}},$\ \par 
because ${\left\Vert{\partial \chi }\right\Vert}_{\infty }\leq R^{-1}.$\ \par 
\quad We can use the "Peter-Paul" inequality~\cite[Theorem 7.28, p.
 173]{GuilbargTrudinger98} (see also~\cite[Theorem 6.18, (g)
 p. 232]{Warner83} for the case $r=2$).\ \par 
\quad \quad \quad $\displaystyle \exists C>0,\ \forall \epsilon >0::{\left\Vert{\chi
 u}\right\Vert}_{W^{m-1,r}({\mathbb{R}}^{n})^{N}}\leq \epsilon
 {\left\Vert{\chi u}\right\Vert}_{W^{m,r}({\mathbb{R}}^{n})^{N}}+C\epsilon
 ^{-m+1}{\left\Vert{\chi u}\right\Vert}_{L^{r}({\mathbb{R}}^{n})^{N}},$\ \par 
with $C$ independent of $R$ of course. We choose $\epsilon =R/2$
 and we get\ \par 
\quad \quad \quad $\displaystyle {\left\Vert{E}\right\Vert}_{L^{r}({\mathbb{R}}^{n})^{N}}\leq
 R^{-1}{\left\Vert{\chi u}\right\Vert}_{W^{m-1,r}({\mathbb{R}}^{n})^{N}}\leq
 \frac{1}{2}{\left\Vert{\chi u}\right\Vert}_{W^{m,r}({\mathbb{R}}^{n})^{N}}+cR^{-m+1}{\left\Vert{\chi
 u}\right\Vert}_{L^{r}({\mathbb{R}}^{n})^{N}}.$\ \par 
Integrating the $s$ power for $t$ in $J$ we get\ \par 
\quad \quad \quad $\displaystyle {\left\Vert{E}\right\Vert}_{L(s,r)}\leq R^{-1}{\left\Vert{\chi
 u}\right\Vert}_{L^{s}(J,W^{m-1,r}({\mathbb{R}}^{n})^{N})}\leq
 \frac{1}{2}{\left\Vert{\chi u}\right\Vert}_{W(s,r)}+cR^{-m+1}{\left\Vert{\chi
 u}\right\Vert}_{L(s,r)}.$\ \par 
Hence putting it in~(\ref{pBB22}), we get:\ \par 
\quad \quad \quad $\displaystyle {\left\Vert{D(\chi u)}\right\Vert}_{L(s,r)}\leq
 {\left\Vert{\chi D(u)}\right\Vert}_{L(s,r)}+\frac{1}{2}{\left\Vert{\chi
 u}\right\Vert}_{W(s,r)}+cR^{-m+1}{\left\Vert{\chi u}\right\Vert}_{L(s,r)}.$\
 \par 
Now using~(\ref{pBB23}) we have, because $R\leq 1\Rightarrow
 R^{-m+1}\leq R^{-m},$\ \par 
\quad \quad \quad $\displaystyle {\left\Vert{\partial _{t}(\chi u)}\right\Vert}_{L(s,r)}+\frac{1}{2}{\left\Vert{\chi
 u}\right\Vert}_{W(s,r)}\leq c_{1}{\left\Vert{\chi Du}\right\Vert}_{L(s,r)}+c_{2}R^{-m}{\left\Vert{\chi
 u}\right\Vert}_{L(s,r)}.$\ \par 
Because $\chi =1$ in $B^{1}$ and $\chi \geq 0$ we get\ \par 
\quad \quad \quad $\displaystyle {\left\Vert{\partial _{t}(u)}\right\Vert}_{L^{s}(J,L^{r}(B^{1})^{N})}+{\left\Vert{u}\right\Vert}_{L^{s}(J,W^{m,r}(B^{1})^{N})}\leq
 {\left\Vert{\partial _{t}(\chi u)}\right\Vert}_{L(s,r)}+{\left\Vert{\chi
 u}\right\Vert}_{W(s,r)}.$\ \par 
And, because $\chi \leq 1$ with compact support in $B,$ we deduce\ \par 
\quad \quad \quad $\displaystyle c_{1}{\left\Vert{\chi Du}\right\Vert}_{L(s,r)}+c_{2}R^{-m}{\left\Vert{\chi
 u}\right\Vert}_{L(s,r)}\leq c_{1}{\left\Vert{Du}\right\Vert}_{L^{s}(J,L^{r}(B)^{N})}+c_{2}R^{-m}{\left\Vert{\chi
 u}\right\Vert}_{L^{s}(J,L^{r}(B)^{N})}.$\ \par 
So finally:\ \par 
\quad \quad \quad $\displaystyle {\left\Vert{\partial _{t}(u)}\right\Vert}_{L^{s}(J,L^{r}(B^{1})^{N})}+{\left\Vert{u}\right\Vert}_{L^{s}(J,W^{m,r}(B^{1})^{N})}\leq
 c_{1}{\left\Vert{Du}\right\Vert}_{L^{s}(J,L^{r}(B)^{N})}+c_{2}R^{-m}{\left\Vert{\chi
 u}\right\Vert}_{L^{s}(J,L^{r}(B)^{N})},$\ \par 
with new constants still not depending on $B$ hence nor on $R.$\ \par 
\quad Up to now we have $J=\lbrack 0,\infty );$ to get a finite interval
 we just multiply $u$ by a function $\psi (t)$ with compact support
 in $\lbrack 0,T+\alpha )$ such that $0\leq \psi \leq 1,\ \psi
 (t)=1$ for $t\in \lbrack 0,T\rbrack $ and, using that\ \par 
\quad \quad \quad $\displaystyle \partial _{t}(\psi u)=\psi 'u+\psi \partial _{t}u\Rightarrow
 {\left\Vert{\partial _{t}(\psi u)}\right\Vert}_{L^{s}(J,L^{r}(B^{1})^{N})}\geq
 {\left\Vert{\psi \partial _{t}u}\right\Vert}_{L^{s}(J,L^{r}(B^{1})^{N})}-{\left\Vert{\psi
 'u}\right\Vert}_{L^{s}(J,L^{r}(B^{1})^{N})}$\ \par 
we get:\ \par 
\quad \quad \quad $\displaystyle {\left\Vert{\psi \partial _{t}u}\right\Vert}_{L^{s}(J,L^{r}(B^{1})^{N})}-{\left\Vert{\psi
 'u}\right\Vert}_{L^{s}(J,L^{r}(B^{1})^{N})}+{\left\Vert{\psi
 u}\right\Vert}_{L^{s}(J,W^{m,r}(B^{1})^{N})}\leq $\ \par 
\quad \quad \quad \quad \quad \quad \quad $\displaystyle \leq c_{1}{\left\Vert{D(\psi u)}\right\Vert}_{L^{s}(J,L^{r}(B)^{N})}+c_{2}R^{-m}{\left\Vert{\psi
 u}\right\Vert}_{L^{s}(J,L^{r}(B)^{N})}.$\ \par 
But, because $\psi $ depends only on $t,\ D(\psi u)=\psi 'u+\psi
 Du$ we have\ \par 
\quad \quad \quad $\displaystyle {\left\Vert{D(\psi u)}\right\Vert}_{L^{s}(J,L^{r}(B)^{N})}={\left\Vert{\psi
 Du}\right\Vert}_{L^{s}(J,L^{r}(B)^{N})}+{\left\Vert{\psi 'u}\right\Vert}_{L^{s}(J,L^{r}(B)^{N})}.$\
 \par 
So we deduce\ \par 
\quad \quad \quad $\displaystyle {\left\Vert{\psi \partial _{t}u}\right\Vert}_{L^{s}(J,L^{r}(B^{1})^{N})}+{\left\Vert{\psi
 u}\right\Vert}_{L^{s}(J,W^{m,r}(B^{1})^{N})}\leq $\ \par 
\quad \quad \quad \quad \quad \quad \quad $\displaystyle \leq c_{1}{\left\Vert{\psi Du}\right\Vert}_{L^{s}(J,L^{r}(B)^{N})}+(1+c_{1}){\left\Vert{\psi
 'u}\right\Vert}_{L^{s}(J,L^{r}(B)^{N})}+c_{2}R^{-m}{\left\Vert{\psi
 u}\right\Vert}_{L^{s}(J,L^{r}(B)^{N})}.$\ \par 
\quad Now we have that   $\left\vert{\psi '}\right\vert \leq C$ and
 $R\leq 1$ so we end with:\ \par 
\quad \quad \quad $\displaystyle {\left\Vert{\partial _{t}u}\right\Vert}_{L^{s}(\lbrack
 0,T\rbrack ,L^{r}(B^{1})^{N})}+{\left\Vert{u}\right\Vert}_{L^{s}(\lbrack
 0,T\rbrack ,W^{m,r}(B^{1})^{N})}\leq $\ \par 
\quad \quad \quad $\displaystyle \leq c_{1}{\left\Vert{Du}\right\Vert}_{L^{s}(\lbrack
 0,T+\alpha \rbrack ,L^{r}(B)^{N})}+c_{2}R^{-m}{\left\Vert{u}\right\Vert}_{L^{s}(\lbrack
 0,T+\alpha \rbrack ,L^{r}(B)^{N})},$\ \par 
the new constants now depend on $\alpha $ (and $\mu $) but still
 not on $B$ hence not on $R.$\ \par 
The proof is complete. $\blacksquare $\ \par 

\subsection{Sobolev comparison estimates.}
\quad The following two lemmas are more or less well known. I give
 the proofs here for the reader convenience.\ \par 

\begin{lem}
~\label{m3}Let $B(x,R)$ be a $(m,\epsilon )$-admissible ball
 in $M$ and $\varphi \ :\ B(x,R)\rightarrow {\mathbb{R}}^{n}$
  be the admissible chart relative to $B(x,R).$ Set $v:=\varphi
 ^{*}u,$ then, for $m\geq 1$:\par 
\quad \quad \quad $\displaystyle \forall u\in W_{G}^{m,r}(B(x,R)),\ {\left\Vert{u}\right\Vert}_{W_{G}^{m,r}(B(x,R))}\leq
 cR^{-m}{\left\Vert{v}\right\Vert}_{W^{m,r}(\varphi (B(x,R)))},$\par 
and, with $B_{e}(0,t)$ the euclidean ball in ${\mathbb{R}}^{n}$
 centered at $0$ and of radius $t,$\par 
\quad \quad \quad $\displaystyle {\left\Vert{v}\right\Vert}_{W^{m,r}(B_{e}(0,(1-\epsilon
 )R))}\leq cR^{-m}{\left\Vert{u}\right\Vert}_{W_{G}^{m,r}(B(x,R))}.$\par 
\quad We also have, for $m=0$:\par 
\quad \quad \quad $\displaystyle \forall u\in L_{G}^{r}(B(x,R)),\ {\left\Vert{u}\right\Vert}_{L_{G}^{r}(B(x,R))}\leq
 (1+C\epsilon ){\left\Vert{v}\right\Vert}_{L^{r}(\varphi (B(x,R)))},$\par 
and\par 
\quad \quad \quad $\displaystyle {\left\Vert{v}\right\Vert}_{L^{r}(B_{e}(0,(1-\epsilon
 )R))}\leq (1+C\epsilon ){\left\Vert{u}\right\Vert}_{L_{G}^{r}(B(x,R))}.$\par 
The constants $c,\ C$ being independent of $B.$\par 
\quad In the case of a function $u$ on $M,$ we have better results.
 Let $B(x,R)$ be a $(m-1,\epsilon )$-admissible ball in $M$ and
 $\varphi \ :\ B(x,R)\rightarrow {\mathbb{R}}^{n}$  be the admissible
 chart relative to $B(x,R).$ Set $v:=u\circ \varphi ,$ then, for $m\geq 1$:\par 
\quad \quad \quad $\displaystyle \forall u\in W^{m,r}(B(x,R)),\ {\left\Vert{u}\right\Vert}_{W^{m,r}(B(x,R))}\leq
 cR^{1-m}{\left\Vert{v}\right\Vert}_{W^{m,r}(\varphi (B(x,R)))},$\par 
and\par 
\quad \quad \quad $\displaystyle {\left\Vert{v}\right\Vert}_{W^{m,r}(B_{e}(0,(1-\epsilon
 )R))}\leq cR^{1-m}{\left\Vert{u}\right\Vert}_{W^{m,r}(B(x,R))}.$\par 
\quad We also have, for $m=0$:\par 
\quad \quad \quad $\displaystyle \forall u\in L^{r}(B(x,R)),\ {\left\Vert{u}\right\Vert}_{L_{G}^{r}(B(x,R))}\leq
 (1+C\epsilon ){\left\Vert{v}\right\Vert}_{L^{r}(\varphi (B(x,R)))},$\par 
and\par 
\quad \quad \quad $\displaystyle {\left\Vert{v}\right\Vert}_{L^{r}(B_{e}(0,(1-\epsilon
 )R))}\leq (1+C\epsilon ){\left\Vert{u}\right\Vert}_{L^{r}(B(x,R))}.$\par 
The constants $c,\ C$ being independent of $B.$
\end{lem}
\quad Proof.\ \par 
We have to compare the norms of $u,\ \nabla u,\cdot \cdot \cdot
 ,\ \nabla ^{m}u,\ $ with the corresponding ones for $v:=\varphi
 ^{*}u$ in ${\mathbb{R}}^{n}.$\ \par 
By Lemma~\ref{SB30} the $(m,\epsilon )$-admissible ball $\displaystyle
 B(x,R)$ trivialises the bundle $G,$ hence the image of a section
 of $G$ in ${\mathbb{R}}^{n}$ is just vectors of functions. Precisely
 $v:=\varphi ^{*}u\in \varphi (B(x,R)){\times}{\mathbb{R}}^{N}.$\ \par 
We have, because $(1-\epsilon )\delta _{ij}\leq g_{ij}\leq (1+\epsilon
 )\delta _{ij}$ in $B(x,R)$:\ \par 
\quad \quad \quad $\displaystyle B_{e}(0,(1-\epsilon )R)\subset \varphi (B(x,R))\subset
 B_{e}(0,(1+\epsilon )R).$\ \par 
Let $u$ be a $G$-form in $M.$ By our assumption (CMT) we have
 that $\nabla u$ depends on the first order derivatives of the
 metric tensor $g.$\ \par 
\quad Because of~(\ref{SB34}) we get, with the fact that $B(x,R)$ is
 $(1,\epsilon )$-admissible, with $\eta :=\frac{\epsilon }{R},$\ \par 
\quad \quad \quad $\displaystyle \sum_{\left\vert{\beta }\right\vert =1}{R\sup
 \ _{i,j=1,...,n,\ y\in B_{x}(R)}\left\vert{\partial ^{\beta
 }g_{ij}(y)}\right\vert }\leq \epsilon \Rightarrow \left\vert{\Gamma
 ^{G,k}_{ij}}\right\vert \leq C\eta ,$\ \par 
with $C$ being independent of $B$ because the constant in the
 (CMT) assumption depends only on $n,\epsilon $ and $G.$\ \par 
Hence\ \par 
\quad \quad \quad $\displaystyle \forall y\in B(x,R),\ \left\vert{u(y)}\right\vert
 =\left\vert{v(z)}\right\vert ,\ \ \left\vert{\nabla u(y)}\right\vert
 \leq \left\vert{\partial u}\right\vert +\left\vert{\Phi }\right\vert ,$\ \par 
where $\Phi $ depends on the coefficients of $u$ and on the first
 order derivatives of the metric tensor $g.$\ \par 
So\ \par 
\quad \quad \quad \begin{equation} \ \left\vert{\nabla u(y)}\right\vert \leq \left\vert{\partial
 v(z)}\right\vert +C\eta \left\vert{v(z)}\right\vert .\label{SC22}\end{equation}\
 \par 
\quad Integrating this we get\ \par 
\quad \quad \quad $\displaystyle \ {\left\Vert{\nabla u(y)}\right\Vert}_{L^{r}(B(x,R))}\leq
 {\left\Vert{\partial v}\right\Vert}_{L^{r}(B_{e}(0,(1+\epsilon
 )R))}+C\eta {\left\Vert{v}\right\Vert}_{L^{r}(B_{e}(0,(1+\epsilon
 )R))}.$\ \par 
The same way for $\nabla ^{k}u$ with $1<k\leq m,$ we have:\ \par 
\quad \quad \quad \begin{equation} \sum_{\left\vert{\beta }\right\vert =k}{\sup
 \ _{i,j=1,...,n,\ y\in B_{x}(R)}\left\vert{\partial ^{\beta
 }g_{ij}(y)}\right\vert }\leq \epsilon /R^{k},\label{bPB58}\end{equation}\ \par 
we get: $\displaystyle \nabla ^{k}u=\partial ^{k}v+\Phi $ and
 by Remark~\ref{pBB24} we have $\Phi $ depends on the derivatives
 of the coefficients of $v$ up to order $k-1$ and on the first
  $k$ order derivatives of the metric tensor $g$ by the (CMT)
 assumption. With ~(\ref{bPB58}) for $j$$\leq k$ we get:\ \par 
\quad $\displaystyle \forall y\in B(x,R),\ \left\vert{\nabla ^{k}u(y)}\right\vert
 \leq \left\vert{\partial ^{k}v(z)}\right\vert +\epsilon (C_{0}R^{-1}\left\vert{v(z)}\right\vert
 +\ C_{1}R^{-2})\left\vert{\partial v(z)}\right\vert +\cdot \cdot
 \cdot +C_{k-1}R^{-k}\left\vert{\partial ^{k-1}v(z)}\right\vert ).$\ \par 
So, using that\ \par 
\quad \quad \quad $\displaystyle {\left\Vert{u}\right\Vert}_{W^{m,r}(B(x,R))}={\left\Vert{\nabla
 ^{m}u}\right\Vert}_{L^{r}(B(x,R))}+\cdot \cdot \cdot +{\left\Vert{\nabla
 u}\right\Vert}_{L^{r}(B(x,R))}+{\left\Vert{u}\right\Vert}_{L^{r}(B(x,R))},$\
 \par 
we get\ \par 
\quad $\displaystyle {\left\Vert{u}\right\Vert}_{W^{m,r}(B(x,R))}\leq
 {\left\Vert{\partial ^{m}v}\right\Vert}_{L^{r}(B_{e}(0,(1+\epsilon
 )R))}+\sum_{k=0}^{m-1}{(1+\epsilon C_{k}R^{-k-1}){\left\Vert{\partial
 ^{k}v}\right\Vert}_{L^{r}(B_{e}(0,(1+\epsilon )R))}}.$\ \par 
\quad Because $R\leq 1$ we get that, for $k\geq 1,$\ \par 
\quad \quad \quad $\displaystyle 1+\epsilon C_{k-1}R^{-k}=R^{-k}(R^{k}+\epsilon
 C_{k-1})\leq R^{-k}(1+\epsilon C_{k-1})\leq cR^{-k},$\ \par 
hence\ \par 
\quad \quad \quad $\displaystyle (1+\epsilon C_{k}R^{-k-1}){\left\Vert{\partial
 ^{k}v}\right\Vert}_{L^{r}(B_{e}(0,(1+\epsilon )R))}\leq cR^{-k-1}{\left\Vert{\partial
 ^{k}v}\right\Vert}_{L^{r}(B_{e}(0,(1+\epsilon )R))},$\ \par 
with a new constant $c$ independent of $B.$ So we get\ \par 
\quad \quad \quad $\displaystyle {\left\Vert{u}\right\Vert}_{W^{m,r}(B(x,R))}\leq
 {\left\Vert{\partial ^{m}v}\right\Vert}_{L^{r}(B_{e}(0,(1+\epsilon
 )R))}+c\sum_{k=0}^{m-1}{R^{-k-1}{\left\Vert{\partial ^{k}v}\right\Vert}_{L^{r}(B_{e}(0,(1+\epsilon
 )R))}},$\ \par 
so\ \par 
\begin{equation} {\left\Vert{u}\right\Vert}_{W^{m,r}(B(x,R))}\leq
 {\left\Vert{\partial ^{m}v}\right\Vert}_{L^{r}(B_{e}(0,(1+\epsilon
 )R))}+cR^{-m}\sum_{k=0}^{m-1}{{\left\Vert{\partial ^{k}v}\right\Vert}_{L^{r}(B_{e}(0,(1+\epsilon
 )R))}},\label{m2}\end{equation}\ \par 
and\ \par 
\quad \quad \quad $\displaystyle {\left\Vert{u}\right\Vert}_{W^{m,r}(B(x,R))}\leq
 cR^{-m}{\left\Vert{v}\right\Vert}_{W^{m,r}(B_{e}(0,(1+\epsilon )R))}.$\ \par 
The same way we get the reverse estimates\ \par 
\quad \quad \quad $\displaystyle {\left\Vert{v}\right\Vert}_{W^{m,r}(B_{e}(0,(1-\epsilon
 )R))}\leq cR^{-m}{\left\Vert{u}\right\Vert}_{W^{m,r}(B(x,R)}.$\ \par 
The case $m=0$ is given by the equation~(\ref{m2}).\ \par 
All the constants here are independent of $B.$\ \par 
\quad In the case of a function $u$ on $M,$ we have a better result.
  In this case we have:  $\displaystyle (\nabla u)_{j}:=\partial
 _{j}u\ $ in local coordinates, so $\displaystyle \left\vert{\nabla
 u(y)}\right\vert \leq \left\vert{\partial v(z)}\right\vert .$\ \par 
While the components of $\nabla ^{2}u$ are given by $\displaystyle
 (\nabla ^{2}u)_{ij}=\partial _{ij}u-\Gamma ^{k}_{ij}\partial
 _{k}u,$ where the Christoffel symbols $\displaystyle \Gamma
 ^{k}_{ij}$ are those of the Levi-Civita connection. Now we have
 for $B(x,R)$ be $(1,\epsilon )$-admissible ball:\ \par 
\quad \quad \quad $\displaystyle \sum_{\left\vert{\beta }\right\vert =1}{R\sup
 \ _{i,j=1,...,n,\ y\in B_{x}(R)}\left\vert{\partial ^{\beta
 }g_{ij}(y)}\right\vert }\leq \epsilon \Rightarrow \left\vert{\Gamma
 ^{k}_{ij}}\right\vert \leq C\epsilon /R,$\ \par 
because, by Lemma~\ref{pBB26}, (CMT) is true for the Levi-Civita
 connection. So we get\ \par 
\quad \quad \quad $\displaystyle \left\vert{\nabla ^{2}u(y)}\right\vert \leq \left\vert{\partial
 ^{2}v(z)}\right\vert +c\frac{\epsilon }{R}\left\vert{\partial
 v(z)}\right\vert .$\ \par 
Hence by integration we get\ \par 
\quad \quad \quad $\displaystyle {\left\Vert{\nabla ^{2}u}\right\Vert}_{L^{r}(B(x,R))}\leq
 {\left\Vert{\partial ^{2}v}\right\Vert}_{L^{r}(B_{e}(0,(1+\epsilon
 )R))}+C\frac{\epsilon }{R}{\left\Vert{v\partial v}\right\Vert}_{L^{r}(B_{e}(0,(1+\epsilon
 )R))}.$\ \par 
For controlling $\nabla ^{k}u$ we need only to have  $B(x,R)$
 be $(k-1,\epsilon )$-admissible and we get:\ \par 
\quad \quad \quad $\displaystyle \forall y\in B(x,R),\ \left\vert{\nabla ^{k}u(y)}\right\vert
 \leq \left\vert{\partial ^{k}v(z)}\right\vert +\epsilon (C_{0}\left\vert{v(z)}\right\vert
 +\ C_{1}R^{-1})\left\vert{\partial v(z)}\right\vert +\cdot \cdot
 \cdot +C_{k-1}R^{1-k}\left\vert{\partial ^{k-1}v(z)}\right\vert ).$\ \par 
Now, the same way as for the sections of $G,$ we get, with $B(x,R)$
 be $(m-1,\epsilon )$-admissible:\ \par 
\quad \quad \quad $\displaystyle {\left\Vert{u}\right\Vert}_{W^{m,r}(B(x,R))}\leq
 cR^{1-m}{\left\Vert{v}\right\Vert}_{W^{m,r}(B_{e}(0,(1+\epsilon )R))}.$\ \par 
We also get the reverse estimates\ \par 
\quad \quad \quad $\displaystyle {\left\Vert{v}\right\Vert}_{W^{m,r}(B_{e}(0,(1-\epsilon
 )R))}\leq cR^{1-m}{\left\Vert{u}\right\Vert}_{W^{m,r}(B(x,R)}.$\ \par 
\quad The proof of the lemma is complete. $\blacksquare $\ \par 

\begin{lem}
~\label{pBB10}(Sobolev embedding) Let $B(x,R)$ is a $(m,\epsilon
 )$-admissible ball in $M$ and $\varphi \ :\ B(x,R)\rightarrow
 {\mathbb{R}}^{n}$  be the admissible chart relative to $B(x,R).$
 We have the Sobolev inequality, for $m\geq 1$:\par 
\quad \quad \quad $\displaystyle \forall u\in W_{G}^{m,\rho }(B(x,R)),\ {\left\Vert{u}\right\Vert}_{L_{G}^{\tau
 }(B(x,R/2))}\leq cR^{-2m}{\left\Vert{u}\right\Vert}_{W_{G}^{m,\rho
 }(B(x,R))}$ with $\displaystyle \frac{1}{\tau }=\frac{1}{\rho
 }-\frac{m}{n}.$\par 
\quad In the special case of functions, with $B(x,R)$ a $(m-1,\epsilon
 )$-admissible ball in $M,$ we have, for $m\geq 1$:\par 
\quad \quad \quad $\displaystyle \forall u\in W^{m,\rho }(B(x,R)),\ {\left\Vert{u}\right\Vert}_{L^{\tau
 }(B(x,R/2))}\leq cR^{1-2m}{\left\Vert{u}\right\Vert}_{W^{m,\rho
 }(B(x,R))}.$\par 
The constant $c$ being independent of $u$ and of the ball $\displaystyle
 B(x,R).$
\end{lem}
\quad Proof.\ \par 
Let us see first what happen in ${\mathbb{R}}^{n}.$ In the ball
 $B_{1}:=B_{e}(0,1)$ we have the Sobolev embedding:\ \par 
\quad \quad \quad $\displaystyle {\left\Vert{u}\right\Vert}_{L^{\tau }(B_{1})}\leq
 C{\left\Vert{u}\right\Vert}_{W^{m,\rho }(B_{1})}.$\ \par 
By the change of variables $y=Rx$ we have $x\in B_{e}(0,1)\Rightarrow
 y\in B_{R}:=B_{e}(0,R)$ and, with $v(y):=u(x),$\ \par 
\quad \quad \quad $\displaystyle \ {\left\Vert{v}\right\Vert}_{L^{\tau }(B_{R})}^{\tau
 }:=\int_{B_{R}}{\left\vert{v(y)}\right\vert ^{\tau }dy}=\int_{B_{1}}{\left\vert{v(Rx)}\right\vert
 ^{\tau }J(x)dx}$\ \par 
where $J(x)=R^{n}$ is the Jacobian. So we get\ \par 
\quad \quad \quad $\displaystyle {\left\Vert{v}\right\Vert}_{L^{\tau }(B_{R})}=R^{n/\tau
 }{\left\Vert{u}\right\Vert}_{L^{\tau }(B_{1})}.$\ \par 
For the derivatives of $v$ we get:\ \par 
\quad \quad \quad $\displaystyle \frac{\partial v}{\partial y_{j}}(y)=\frac{\partial
 u(x)}{\partial x_{j}}{\times}\frac{\partial x_{j}}{\partial
 y_{j}}=R^{-1}\frac{\partial u(x)}{\partial x_{j}}.$\ \par 
For the derivatives of order $k$ we get, with $\left\vert{\alpha
 }\right\vert =k,\ \partial ^{\alpha }v=R^{-k}\partial ^{\alpha }u.$\ \par 
Using that\ \par 
\quad \quad \quad $\displaystyle {\left\Vert{u}\right\Vert}_{W^{m,\rho }(B_{1})}=\sum_{k=0,...,m;\left\vert{\alpha
 }\right\vert =k}{{\left\Vert{\partial ^{\alpha }u}\right\Vert}_{L^{\rho
 }(B_{1})}}$\ \par 
and the analogous for $\displaystyle {\left\Vert{v}\right\Vert}_{W^{m,\rho
 }(B_{R})},$ we get\ \par 
\quad \quad \quad $\displaystyle {\left\Vert{v}\right\Vert}_{W^{m,\rho }(B_{R})}=R^{n/\rho
 }\sum_{k=0,...,m;\left\vert{\alpha }\right\vert =k}{R^{-k}{\left\Vert{\partial
 ^{\alpha }u}\right\Vert}_{L^{\rho }(B_{1})}}.$\ \par 
So we get, with $R\leq 1,$\ \par 
\quad \quad \quad $\displaystyle R^{n/\rho }{\left\Vert{u}\right\Vert}_{W^{m,\rho
 }(B_{1})}\leq {\left\Vert{v}\right\Vert}_{W^{m,\rho }(B_{R})}\leq
 R^{n/\rho -m}{\left\Vert{u}\right\Vert}_{W^{m,\rho }(B_{1})}.$\ \par 
Hence, from $\displaystyle {\left\Vert{u}\right\Vert}_{L^{\tau
 }(B_{1})}\leq C{\left\Vert{u}\right\Vert}_{W^{m,\rho }(B_{1})}$ we get:\ \par 
\quad \quad \quad $\displaystyle {\left\Vert{v}\right\Vert}_{L^{\tau }(B_{R})}=R^{n/\tau
 }{\left\Vert{u}\right\Vert}_{L^{\tau }(B_{1})}\leq CR^{n/\tau
 }{\left\Vert{u}\right\Vert}_{W^{m,\rho }(B_{1})}\leq $\ \par 
\quad \quad \quad \quad \quad \quad \quad $\displaystyle \leq CR^{n/\tau -n/\rho }R^{n/\rho }{\left\Vert{u}\right\Vert}_{W^{m,\rho
 }(B_{1})}\leq CR^{-m}{\left\Vert{v}\right\Vert}_{W^{m,\rho }(B_{R})},$\ \par 
because $\displaystyle \frac{1}{\tau }=\frac{1}{\rho }-\frac{m}{n}.$\ \par 
I.e. with $R\leq 1,$\ \par 
\quad \quad \quad \begin{equation} {\left\Vert{v}\right\Vert}_{L^{\tau }(B_{R})}\leq
 CR^{-m}{\left\Vert{v}\right\Vert}_{W^{m,\rho }(B_{R})}.\label{pBB9}\end{equation}\
 \par 
\quad We shall use Lemma~\ref{m3} to get the precise Sobolev embeddings
 we need in the manifold $M.$\ \par 
Let $B(x,R)$ be a $(m,\epsilon )$-admissible ball in $M$ and
 $\varphi \ :\ B(x,R)\rightarrow {\mathbb{R}}^{n}$  be the admissible
 chart relative to $B(x,R).$ Set $v:=\varphi ^{*}u,$ then, for
 $m\geq 1$ the Lemma~\ref{m3} gives:\ \par 
\quad \quad \quad $\displaystyle \forall u\in L_{G}^{\tau }(B(x,R)),\ {\left\Vert{u}\right\Vert}_{L_{G}^{\tau
 }(B(x,R))}\leq (1+C\epsilon ){\left\Vert{v}\right\Vert}_{L^{\tau
 }(\varphi (B(x,R)))}\leq c{\left\Vert{v}\right\Vert}_{L^{\tau
 }(B_{e}(0,(1+\epsilon )R))}.$\ \par 
And also:\ \par 
\quad \quad \quad $\displaystyle {\left\Vert{v}\right\Vert}_{W^{m,\rho }(B_{e}(0,(1-\epsilon
 )R))}\leq cR^{-m}{\left\Vert{u}\right\Vert}_{W_{G}^{m,\rho }(B(x,R))}.$\ \par 
Putting this in~(\ref{pBB9}) we get,\ \par 
\quad \quad \quad $\displaystyle {\left\Vert{u}\right\Vert}_{L_{G}^{\tau }(B(x,R/2))}\leq
 c{\left\Vert{v}\right\Vert}_{L^{\tau }(B_{e}(0,(1+\epsilon )R/2))}$\ \par 
and\ \par 
\quad \quad \quad $\displaystyle {\left\Vert{v}\right\Vert}_{L^{\tau }(B_{e}(0,(1+\epsilon
 )R/2))}\leq CR^{-m}{\left\Vert{v}\right\Vert}_{W^{m,\rho }(B_{e}(0,(1-\epsilon
 )R))}\leq cR^{-2m}{\left\Vert{u}\right\Vert}_{W_{G}^{m,\rho }(B(x,R))}.$\ \par 
provided that $(1+\epsilon )\frac{R}{2}\leq (1-\epsilon )R\iff
 \epsilon \leq \frac{1}{3}.$ Then:\ \par 
\quad \quad \quad $\displaystyle {\left\Vert{u}\right\Vert}_{L_{G}^{\tau }(B(x,R/2))}\leq
 cR^{-2m}{\left\Vert{u}\right\Vert}_{W_{G}^{m,\rho }(B(x,R))}$
 with $\displaystyle \frac{1}{\tau }=\frac{1}{\rho }-\frac{m}{n}.$\ \par 
\quad In the special case of functions we need only to have $B(x,R)$
 be a $(m-1,\epsilon )$-admissible ball in $M,$ and we get:\ \par 
\quad \quad \quad $\displaystyle \forall u\in L^{\tau }(B(x,R)),\ {\left\Vert{u}\right\Vert}_{L^{\tau
 }(B(x,R))}\leq c{\left\Vert{v}\right\Vert}_{L^{\tau }(B_{e}(0,(1+\epsilon
 )R))}$\ \par 
and\ \par 
\quad \quad \quad $\displaystyle {\left\Vert{u}\right\Vert}_{L^{\tau }(B(x,R/2))}\leq
 cR^{1-2m}{\left\Vert{u}\right\Vert}_{W^{m,\rho }(B(x,R))}.$\ \par 
The proof is complete. $\blacksquare $\ \par 

\subsection{The main local estimates.}
\quad We shall use the following notation to ease the writing:\ \par 

\begin{defin}
~\label{pBB13} For $r,s>1,\ \alpha >0$ fixed and $m,k\in {\mathbb{N}},\
 m\geq 2,$ we set:\par 
\quad \quad \quad $\displaystyle L(r,k):=L^{s}(\lbrack 0,T+\alpha /2^{k}\rbrack
 ,L_{G}^{r}(B^{k}))$  and $\displaystyle W(r,k):=L^{s}(\lbrack
 0,T+\alpha /2^{k}\rbrack ,W_{G}^{m,r}(B^{k})),$\par 
where $B:=B(x,R)$ is a ball in the riemann manifold $(M,g)$ and
 $B^{k}:=B(x,R/2^{k}).$
\end{defin}
\ \par 
\quad The following theorem follows by standard techniques but is needed
 for the sequel.\ \par 

\begin{thm}
~\label{m6}Suppose that $A$ is a $(C,\theta )$-elliptic operator
 of order $m$ acting on sections of the vector bundle $\displaystyle
 G:=(H,\pi ,M)$ in the complete riemannian manifold $(M,g),$
 with $\theta <\pi /2,$ and consider the parabolic equation $\displaystyle
 Du=\partial _{t}u-Au$ also acting on sections of $G$ and verifying
 $Du\in L^{s}(\lbrack 0,T+\alpha \rbrack ,L_{G}^{r}(B))$ and
 $u\in L^{s}(\lbrack 0,T+\alpha \rbrack ,L_{G}^{r}(B)).$\par 
Let $\displaystyle B:=B(x,R)$ be a $(m,\epsilon )$-admissible
 ball and set $\displaystyle B^{1}:=B(x,R/2).$ Then, with $r,s$
 in $(1,\infty ),$ we have:\par 
\quad \quad \quad $\displaystyle {\left\Vert{\partial _{t}u}\right\Vert}_{L^{s}(\lbrack
 0,T+\alpha /2\rbrack ,L_{G}^{r}(B^{1}))}+R^{m}{\left\Vert{u}\right\Vert}_{L^{s}(\lbrack
 0,T+\alpha /2\rbrack ,W_{G}^{m,r}(B^{1}))}\leq $\par 
\quad \quad \quad \quad \quad \quad \quad $\displaystyle \leq c_{3}{\left\Vert{Du}\right\Vert}_{L^{s}(\lbrack
 0,T+\alpha \rbrack ,L_{G}^{r}(B))}+c_{4}R^{-m}{\left\Vert{u}\right\Vert}_{L^{s}(\lbrack
 0,T+\alpha \rbrack ,L_{G}^{r}(B))}.$\par 
In the case of functions we get, with this time $\displaystyle
 B\in {\mathcal{A}}_{m-1}(\epsilon ),$\par 
\quad \quad \quad $\displaystyle {\left\Vert{\partial _{t}u}\right\Vert}_{L^{s}(\lbrack
 0,T+\alpha /2\rbrack ,L_{G}^{r}(B^{1}))}+R^{m-1}{\left\Vert{u}\right\Vert}_{L^{s}(\lbrack
 0,T+\alpha /2\rbrack ,W_{G}^{m,r}(B^{1}))}\leq $\par 
\quad \quad \quad \quad \quad $\displaystyle \leq c_{3}{\left\Vert{Du}\right\Vert}_{L^{s}(\lbrack
 0,T+\alpha \rbrack ,L_{G}^{r}(B))}+c_{4}R^{-m}{\left\Vert{u}\right\Vert}_{L^{s}(\lbrack
 0,T+\alpha \rbrack ,L_{G}^{r}(B))}.$\par 
The constants $c_{3},\ c_{4}$ are independent of $u$ and of $B.$
\end{thm}
\quad Proof.\ \par 
The ball $B$ being admissible, there is a diffeomorphism $\varphi
 :B\rightarrow {\mathbb{R}}^{n}$ such that $G$ trivialises on
 $B.$ I.e. we have, for any section $u$ over $B$:\ \par 
\quad \quad \quad $\displaystyle \pi ^{-1}(B)\rightarrow B{\times}H,\ u\rightarrow
 (\pi (u),\chi _{\varphi }(u)).$\ \par 
So the local representation of the section $u$ is: $\displaystyle
 u_{\varphi }:=\chi _{\varphi }\circ u\circ \varphi ^{-1}.$\ \par 
We shall apply Theorem~\ref{m1} with a slight change in $T$ and
 $\alpha $ to the images of $A,G,u,$\ \par 
\quad $\displaystyle (*)\ \ \ \ \ \ \ \ {\left\Vert{\partial _{t}u_{\varphi
 }}\right\Vert}_{L(r,1)}+{\left\Vert{u_{\varphi }}\right\Vert}_{W(r1)}\leq
 c_{1}{\left\Vert{(Du)_{\varphi }}\right\Vert}_{L(r,0)}+c_{2}R_{\varphi
 }^{-m}{\left\Vert{u_{\varphi }}\right\Vert}_{L(r,0)},$\ \par 
where $A_{\varphi },B_{\varphi },R_{\varphi },u_{\varphi }$ are
 the images by $\varphi $ of $A,B,R,u$ and the image of $G$ is
 the trivial bundle  $\varphi (B){\times}{\mathbb{R}}^{N}$ in
 $\displaystyle {\mathbb{R}}^{n}.$ The constants $c_{1},\ c_{2}$
 being independent of $B_{\varphi }.$\ \par 
\quad First, because of the condition $\displaystyle (1-\epsilon )\delta
 _{ij}\leq g_{ij}\leq (1+\epsilon )\delta _{ij}$ in the definition
 of the $\epsilon $-admissible ball, we have that $R_{\varphi }\simeq R.$\ \par 
\quad Now we use the Sobolev comparison estimates given by Lemma~\ref{m3}
 to get:\ \par 
\quad \quad \quad $\displaystyle {\left\Vert{\partial _{t}u}\right\Vert}_{L_{G}^{r}(B^{1})}\leq
 (1+C\epsilon ){\left\Vert{\partial _{t}u_{\varphi }}\right\Vert}_{L^{r}(\varphi
 (B^{1}))},$\ \par 
because $(\partial _{t}u)_{\varphi }=\partial _{t}u_{\varphi
 }.$ We also have:\ \par 
\quad \quad \quad $\displaystyle R^{m}{\left\Vert{u}\right\Vert}_{W_{G}^{m,r}(B^{1})}\leq
 c{\left\Vert{u_{\varphi }}\right\Vert}_{W^{m,r}(\varphi (B^{1}))}.$\ \par 
The constants $c,\ C$ being independent of $B.$\ \par 
Integrating the $s$-power with respect to $t,$ we get for the
 left hand side of (*)\ \par 
\quad \quad \quad $\displaystyle {\left\Vert{\partial _{t}u}\right\Vert}_{L(r,1)}+R^{m}{\left\Vert{u}\right\Vert}_{W(r,1)}\leq
 C({\left\Vert{\partial _{t}u_{\varphi }}\right\Vert}_{L(r,1)}+{\left\Vert{u_{\varphi
 }}\right\Vert}_{W(r,1)}).$\ \par 
\quad Now, still by Lemma~\ref{m3}\ \par 
\quad \quad \quad $\displaystyle {\left\Vert{(Du)_{\varphi }}\right\Vert}_{L^{r}(B_{\varphi
 })^{N}}\leq (1+C\epsilon ){\left\Vert{Du}\right\Vert}_{L_{G}^{r}(B)}$\ \par 
and\ \par 
\quad \quad \quad $\displaystyle {\left\Vert{u_{\varphi }}\right\Vert}_{L^{r}(B_{\varphi
 })^{N}}\leq (1+C\epsilon ){\left\Vert{u}\right\Vert}_{L_{G}^{r}(B)}.$\ \par 
Again integrating the $s$-power with respect to $t,$ we get for
 the right hand side of (*)\ \par 
\quad \quad \quad $\displaystyle c_{1}{\left\Vert{(Du)_{\varphi }}\right\Vert}_{L(r,0)}+c_{2}R_{\varphi
 }^{-m}{\left\Vert{u_{\varphi }}\right\Vert}_{L(r,0)}\leq c_{3}{\left\Vert{Du}\right\Vert}_{L(r,0)}+c_{4}R^{-m}{\left\Vert{u}\right\Vert}_{L(r,0)}.$\
 \par 
Hence replacing in (*) we get, with new constants:\ \par 
\quad \quad \quad $\displaystyle {\left\Vert{\partial _{t}u}\right\Vert}_{L(r,1)}+R^{m}{\left\Vert{u}\right\Vert}_{W(r,1)}\leq
 c_{3}{\left\Vert{Du}\right\Vert}_{L(r,0)}+c_{4}R^{-m}{\left\Vert{u}\right\Vert}_{L(r,0)}.$\
 \par 
The constants $c_{3},\ c_{4}$ are still independent of $u$ and
 of $B$ but depend on $T$ and $\alpha .$\ \par 
\ \par 
\quad In the case of functions, using Lemma~\ref{m3}, we get with this
 time $\displaystyle B\in {\mathcal{A}}_{m-1}(\epsilon ),$\ \par 
\quad \quad \quad $\displaystyle {\left\Vert{\partial _{t}u}\right\Vert}_{L(r,1)}+R^{m-1}{\left\Vert{u}\right\Vert}_{W(r,1)}\leq
 c_{3}{\left\Vert{Du}\right\Vert}_{L(r,0)}+c_{4}R^{-m}{\left\Vert{u}\right\Vert}_{L(r,0)}.$\
 \par 
The proof is complete. $\blacksquare $\ \par 
\ \par 
\quad The following corollary, the LIR inequality, is at the heart
 of the method we use. The induction step works because of the
 gain in regularity we get by this corollary.\ \par 

\begin{cor}
~\label{pBB11}(The LIR inequality) Suppose that $A$ is a $(C,\theta
 )$-elliptic operator of order $m$ acting on sections of  the
 adapted vector bundle $\displaystyle G:=(H,\pi ,M)$ in the complete
 riemannian manifold $(M,g),$ with $\theta <\pi /2,$ and consider
 the parabolic equation $\displaystyle Du=\partial _{t}u-Au$
 also acting on sections of $G.$\par 
Let $\displaystyle B:=B(x,R)$ be a $(m,\epsilon )$-admissible
 ball and set $\displaystyle B^{k}:=B(x,R/2^{k}).$ Then, with
 $r,s$ in $(1,\infty ),$ and $\alpha >0,$ we have:\par 
\quad $\displaystyle R^{m}{\left\Vert{\partial _{t}u}\right\Vert}_{L^{s}(\lbrack
 0,T+\alpha /2^{k+1}\rbrack ,L_{G}^{r}(B^{k+1}))}+R^{2m}{\left\Vert{u}\right\Vert}_{L^{s}(\lbrack
 0,T+\alpha /2^{k+1}\rbrack ,W_{G}^{m,r}(B^{k+1}))}\leq $\par 
\quad \quad \quad $\displaystyle \leq c_{3}R^{m}{\left\Vert{Du}\right\Vert}_{L^{s}(\lbrack
 0,T+\alpha /2^{k}\rbrack ,L_{G}^{r}(B^{k}))}+c_{4}{\left\Vert{u}\right\Vert}_{L^{s}(\lbrack
 0,T+\alpha /2^{k}\rbrack ,L_{G}^{r}(B^{k}))}.$\par 
With the notation of Definition~\ref{pBB13} this gives:\par 
\quad \quad \quad $\displaystyle R^{m}{\left\Vert{\partial _{t}u}\right\Vert}_{L(r,k+1)}+R^{2m}{\left\Vert{u}\right\Vert}_{W(r,k+1)}\leq
 c_{3}R^{m}{\left\Vert{Du}\right\Vert}_{L(r,k)}+c_{4}{\left\Vert{u}\right\Vert}_{L(r,k)}.$\par
 
In the case of functions instead of sections of $G,$ we have,
 with $B(x,R)\in {\mathcal{A}}_{m-1}(\epsilon ),$\par 
\quad $\displaystyle R^{m}{\left\Vert{\partial _{t}u}\right\Vert}_{L^{s}(\lbrack
 0,T+\alpha /2^{k+1}\rbrack ,L^{r}(B^{k+1}))}+R^{2m-1}{\left\Vert{u}\right\Vert}_{L^{s}(\lbrack
 0,T+\alpha /2\rbrack ,W^{m,r}(B^{k+1}))}\leq $\par 
\quad \quad \quad $\displaystyle \leq c_{3}R^{m}{\left\Vert{Du}\right\Vert}_{L^{s}(\lbrack
 0,T+\alpha /2^{k}\rbrack ,L^{r}(B^{k}))}+c_{4}{\left\Vert{u}\right\Vert}_{L^{s}(\lbrack
 0,T+\alpha /2^{k}\rbrack ,L^{r}(B^{k}))}.$\par 
The constants $c_{3},c_{4}$ being independent of $u$ and $B,$
 but depend on $T$ and $\alpha ,$ hence on $k.$
\end{cor}
\quad Proof.\ \par 
We apply Theorem~\ref{m6} to $B^{k+1}\subset B^{k}$ instead of
 $B^{1}\subset B$ and with $\alpha /2^{k}$ instead of $\alpha
 .$ $\blacksquare $\ \par 

\subsection{The induction.}

\begin{rem}
The idea under this method is the following one.\par 
If we have a $\displaystyle u\in L^{s}(\lbrack 0,T\rbrack ,L^{\rho
 }_{G}(B))::Du=\omega $ then the LIR, Corollary~\ref{pBB11},
 gives essentially that $\displaystyle u\in L^{s}(\lbrack 0,T\rbrack
 ,W^{m,\rho }_{G}(B)).$ By applying the Sobolev embedding, Lemma~\ref{pBB10},
 we get $\displaystyle u\in L^{s}(\lbrack 0,T\rbrack ,L^{\tau
 }_{G}(B))::Du=\omega ,$ with $\displaystyle \frac{1}{\tau }=\frac{1}{\rho
 }-\frac{m}{n}.$ But if $\displaystyle \omega \in L^{s}(\lbrack
 0,T\rbrack ,L^{\tau }_{G}(B))$ then a new application of the
 LIR gives $\displaystyle u\in L^{s}(\lbrack 0,T\rbrack ,W^{m,\tau
 }_{G}(B)).$ So we have a strict increase of the regularity of
 $u.$ We can repeat the process up to reach the best regularity
 of the data $\omega .$
\end{rem}
\quad The following lemma is essentially computational.\ \par 

\begin{lem}
(Induction) ~\label{pL30}Provided that:\par 
$\displaystyle IH(k)\ \ \ \ \ \ \ \ \ \ \ \ R^{d_{k}}{\left\Vert{\partial
 _{t}u}\right\Vert}_{L(r_{k},k)}+R^{b_{k}}{\left\Vert{u}\right\Vert}_{W(r_{k},k)}\leq
 c_{1}(k)R^{a_{k}}{\left\Vert{\omega }\right\Vert}_{L(r,0)}+c_{2}(k){\left\Vert{u}\right\Vert}_{L(2,0)}.$\par
 
We get\par 
$\displaystyle IH(k+1)\ \ \ \ \ R^{d_{k+1}}{\left\Vert{\partial
 _{t}u}\right\Vert}_{L(\tau ,k+2)}+R^{b_{k+1}}{\left\Vert{u}\right\Vert}_{W(\tau
 ,k+2)}\leq c_{1}(k+1)R^{a_{k+1}}{\left\Vert{\omega }\right\Vert}_{L(r,0)}+c_{2}(k+1){\left\Vert{u}\right\Vert}_{L(2,0)}.$\par
 
with $\frac{1}{r_{k+1}}=\frac{1}{r_{k}}-\frac{m}{n}=\frac{1}{2}-(k+1)\frac{m}{n},\
 \tau :=\min (r_{k+1},r),$\par 
and for sections of $G$ with $\displaystyle B\in {\mathcal{A}}_{m}(\epsilon
 ),$\par 
\quad \quad \quad $\displaystyle d_{k+1}=3m+b_{k};\ b_{k+1}=4m+b_{k};\ a_{k+1}=\min
 (a_{k},3m+b_{k}),$\par 
and\par 
\quad \quad \quad $\displaystyle c_{1}(k+1)=c_{3}(k)+cc_{4}(k)c_{1}(k);\ c_{2}(k+1)=cc_{4}(k)c_{2}(k).$\par
 
And for functions with $B\in {\mathcal{A}}_{m-1}(\epsilon ),$\par 
\quad \quad \quad $\displaystyle d_{k+1}=3m-1+b_{k};\ b_{k+1}=4m-1+b_{k};\ a_{k+1}=\min
 (a_{k},3m-1+b_{k}),$\par 
and\par 
\quad \quad \quad $\displaystyle c_{1}(k+1)=c_{3}(k)+cc_{4}(k)c_{1}(k);\ c_{2}(k+1)=cc_{4}(k)c_{2}(k).$
\end{lem}

      Proof.\ \par 
We have, by the Sobolev embedding, Lemma~\ref{pBB10}, with $\tau
 :=r_{k+1},\ \rho :=r_{k}$ and $\frac{1}{r_{k+1}}=\frac{1}{r_{k}}-\frac{m}{n},$\
 \par 
\quad \quad \quad $\displaystyle {\left\Vert{u(t,\cdot )}\right\Vert}_{L_{G}^{r_{k+1}}(B^{k+1})}\leq
 cR^{-2m}{\left\Vert{u(t,\cdot )}\right\Vert}_{W_{G}^{m,r_{k}}(B^{k})}$\ \par 
hence, integrating,\ \par 
\quad \quad \quad $\displaystyle {\left\Vert{u}\right\Vert}_{L^{s}(\lbrack 0,T+\alpha
 /2^{k+1}\rbrack ,L^{r_{k+1}}_{p}(B^{k+1}))}\leq cR^{-2m}{\left\Vert{u}\right\Vert}_{L^{s}(\lbrack
 0,T+\alpha /2^{k}\rbrack ,W^{m,r_{k}}_{p}(B^{k}))}.$\ \par 
With the notation of Definition~\ref{pBB13} this gives:\ \par 
\quad \quad \quad \begin{equation} {\left\Vert{u}\right\Vert}_{L(r_{k+1},k+1)}\leq
 cR^{-2m}{\left\Vert{u}\right\Vert}_{W(r_{k},k)}.\label{pBB12}\end{equation}\
 \par 
But by $IH(k)$\ \par 
\quad \quad \quad $\displaystyle R^{b_{k}}{\left\Vert{u}\right\Vert}_{W(r_{k},k)}\leq
 c_{1}(k)R^{a_{k}}{\left\Vert{\omega }\right\Vert}_{L(r,0)}+c_{2}(k){\left\Vert{u}\right\Vert}_{L(2,0)}.$\
 \par 
so\ \par 
\quad $\displaystyle {\left\Vert{u}\right\Vert}_{L(r_{k+1},k+1)}\leq
 cR^{-2m}{\left\Vert{u}\right\Vert}_{W(r_{k},k)}\leq cc_{1}(k)R^{-2m+a_{k}-b_{k}}{\left\Vert{\omega
 }\right\Vert}_{L(r,0)}+cc_{2}(k)R^{-2m-b_{k}}{\left\Vert{u}\right\Vert}_{L(2,0)}.$\
 \par 
\ \par 
Now the LIR inequality, Corollary~\ref{pBB11}, with $\tau =\min
 (r,r^{k+1}),$ gives:\ \par 
\quad \quad \quad $\displaystyle R^{m}{\left\Vert{\partial _{t}u}\right\Vert}_{L(\tau
 ,k+2)}+R^{2m}{\left\Vert{u}\right\Vert}_{W(\tau ,k+2)}\leq c_{3}R^{m}{\left\Vert{Du}\right\Vert}_{L(\tau
 ,k+1)}+c_{4}{\left\Vert{u}\right\Vert}_{L(\tau ,k+1)}.$\ \par 
hence\ \par 
\quad \quad \quad $\displaystyle R^{m}{\left\Vert{\partial _{t}u}\right\Vert}_{L(\tau
 ,k+2)}+R^{2m}{\left\Vert{u}\right\Vert}_{W(\tau ,k+2)}\leq c_{3}R^{m}{\left\Vert{Du}\right\Vert}_{L(\tau
 ,k+1)}+$\ \par 
\quad \quad \quad \quad \quad $\displaystyle +c_{4}cc_{1}(k)R^{-2m+a_{k}-b_{k}}{\left\Vert{\omega
 }\right\Vert}_{L(r,0)}+c_{4}cc_{2}(k)R^{-2m-b_{k}}{\left\Vert{u}\right\Vert}_{L(2,0)}$\
 \par 
because ${\left\Vert{u}\right\Vert}_{L(\tau ,k+1)}\leq {\left\Vert{u}\right\Vert}_{L(r_{k+1},k+1)}.$\
 \par 
But $\tau \leq r,\ \lbrack 0,T+\alpha /2^{k+1}\rbrack \subset
 \lbrack 0,T+\alpha \rbrack ,\ B^{k+2}\subset B,$ so we get\ \par 
\quad \quad \quad $\displaystyle {\left\Vert{Du}\right\Vert}_{L(\tau ,k+1)}\leq
 {\left\Vert{Du}\right\Vert}_{L(r,0)}={\left\Vert{\omega }\right\Vert}_{L(r,0)}.$\
 \par 
Hence\ \par 
\quad \quad \quad $\displaystyle R^{m}{\left\Vert{\partial _{t}u}\right\Vert}_{L(\tau
 ,k+2)}+R^{2m}{\left\Vert{u}\right\Vert}_{W(\tau ,k+2)}\leq $\ \par 
\quad \quad \quad \quad \quad $\displaystyle \leq (c_{3}R^{m}+c_{4}cc_{1}(k))R^{-2m+a_{k}-b_{k}}{\left\Vert{\omega
 }\right\Vert}_{L(r,0)}+c_{4}cc_{2}(k)R^{-2m-b_{k}}{\left\Vert{u}\right\Vert}_{L(2,0)}$\
 \par 
So, multiplying by $R^{2m+b_{k}},$ we get:\ \par 
\quad \quad \quad $\displaystyle R^{3m+b_{k}}{\left\Vert{\partial _{t}u}\right\Vert}_{L(\tau
 ,k+2)}+R^{4m+b_{k}}{\left\Vert{u}\right\Vert}_{W(\tau ,k+2)}\leq $\ \par 
\quad \quad \quad \quad \quad $\displaystyle \leq (c_{3}R^{3m+b_{k}}+c_{4}cc_{1}(k)R^{a_{k}}){\left\Vert{\omega
 }\right\Vert}_{L(r,0)}+c_{4}cc_{2}(k){\left\Vert{u}\right\Vert}_{L(2,0)}$\
 \par 
Hence with\ \par 
\quad \quad \quad $\displaystyle d_{k+1}=3m+b_{k};\ b_{k+1}=4m+b_{k};\ a_{k+1}=\min
 (a_{k},3m+b_{k}),$\ \par 
and\ \par 
\quad \quad \quad $\displaystyle c_{1}(k+1)=c_{3}(\ k)+cc_{4}(k)c_{1}(k);\ c_{2}(k+1)=cc_{4}(k)c_{2}(k),$\
 \par 
we get\ \par 
\quad $\displaystyle IH(k+1)\ \ \ \ \ R^{d_{k+1}}{\left\Vert{\partial
 _{t}u}\right\Vert}_{L(\tau ,k+2)}+R^{b_{k+1}}{\left\Vert{u}\right\Vert}_{W(\tau
 ,k+2)}\leq c_{1}(k+1)R^{a_{k+1}}{\left\Vert{\omega }\right\Vert}_{L(r,0)}+c_{2}(k+1){\left\Vert{u}\right\Vert}_{L(2,0)}.$\
 \par 
\ \par 
\quad In the case of functions, applying again Corollary~\ref{pBB11}, with\ \par 
\quad \quad \quad $\displaystyle d_{k+1}=2m-1+b_{k};\ b_{k+1}=3m-2+b_{k};\ a_{k+1}=\min
 (a_{k},2m-1+b_{k}),$\ \par 
and\ \par 
\quad \quad \quad $\displaystyle c_{1}(k+1)=c_{3}(k)+cc_{4}(k)c_{1}(k);\ c_{2}(k+1)=cc_{4}(k)c_{2}(k),$\
 \par 
we get\ \par 
\quad $\displaystyle IH(k+1)\ \ \ \ \ R^{d_{k+1}}{\left\Vert{\partial
 _{t}u}\right\Vert}_{L(\tau ,k+2)}+R^{b_{k+1}}{\left\Vert{u}\right\Vert}_{W(\tau
 ,k+2)}\leq c_{1}(k+1)R^{a_{k+1}}{\left\Vert{\omega }\right\Vert}_{L(r,0)}+c_{2}(k+1){\left\Vert{u}\right\Vert}_{L(2,0)}.$\
 \par 
The proof is complete. $\blacksquare $\ \par 

\begin{lem}
~\label{pBB14}Let $B=B(x,R)$ be a $\epsilon $-admissible ball
 in $M.$ We have, for $\displaystyle \omega \in L^{s}(\lbrack
 0,T+\alpha \rbrack ,L^{r}_{G}(B))$ with $r\geq 2$:\par 
\quad \quad \quad $\displaystyle {\left\Vert{\omega }\right\Vert}_{L_{G}^{s}(\lbrack
 0,T+\alpha \rbrack ,L^{2}(B))}\leq c(n,\epsilon )R^{\frac{n}{2}-\frac{n}{r}}{\left\Vert{\omega
 }\right\Vert}_{L_{G}^{s}(\lbrack 0,T+\alpha \rbrack ,L^{r}(B))},$\par 
with $c$ depending only on $n$ and $\epsilon .$
\end{lem}
\quad Proof.\ \par 
Let $\omega \in L^{s}(\lbrack 0,T+\alpha \rbrack ,L^{r}_{G}(B)).$
 Because $r\geq 2$ and $B$ is relatively compact, we have $\omega
 \in L^{s}(\lbrack 0,T+\alpha \rbrack ,L^{2}_{G}(B)).$ Because
 $\frac{dv}{\left\vert{B}\right\vert }$ is a probability measure
 on $B,$ where $\left\vert{B}\right\vert $ is the volume of the
 ball $B,$ we get\ \par 
\quad \quad \quad $\displaystyle {\left({\int_{B}{\left\vert{\omega (t,y)}\right\vert
 ^{2}\frac{dv(y)}{\left\vert{B}\right\vert }}}\right)}^{1/2}\leq
 {\left({\int_{B}{\left\vert{\omega (t,y)}\right\vert ^{r}\frac{dv(y)}{\left\vert{B}\right\vert
 }}}\right)}^{1/r},$\ \par 
hence\ \par 
\quad \quad \quad $\displaystyle {\left\Vert{\omega (t,\cdot )}\right\Vert}_{L^{2}(B)}\leq
 \left\vert{B}\right\vert ^{\frac{1}{2}-\frac{1}{r}}{\left\Vert{\omega
 (t,\cdot )}\right\Vert}_{L^{r}(B)}.$\ \par 
Integrating on $t,$ we get\ \par 
\quad \quad \quad $\displaystyle {\left\Vert{\omega }\right\Vert}_{L^{s}(\lbrack
 0,T+\alpha \rbrack ,L^{2}(B))}\leq \left\vert{B}\right\vert
 ^{\frac{1}{2}-\frac{1}{r}}{\left\Vert{\omega }\right\Vert}_{L^{s}(\lbrack
 0,T+\alpha \rbrack ,L^{r}(B))}.$\ \par 
Now on the manifold $M,$ for $B_{x}:=B(x,R)$ a $\epsilon $-admissible
 ball, we get\ \par 
\quad \quad \quad $\displaystyle \forall y\in B_{x},\ (1-\epsilon )^{n}\leq \left\vert{\mathrm{d}\mathrm{e}\mathrm{t}g(y)}\right\vert
 \leq (1+\epsilon )^{n},$\ \par 
hence we have, comparing the Lebesgue measure in ${\mathbb{R}}^{n}$
 with the volume measure in $M,$\ \par 
\quad \quad \quad $\displaystyle \forall x\in M,\ (1-\epsilon )^{n/2}\nu _{n}R^{n}\leq
 \mathrm{V}\mathrm{o}\mathrm{l}(B(x,\ R_{\epsilon }(x)))\leq
 (1+\epsilon )^{n/2}\nu _{n}R^{n},$\ \par 
so, on the manifold $M,$ we have\ \par 
\quad \quad \quad $\displaystyle {\left\Vert{\omega }\right\Vert}_{L^{s}(\lbrack
 0,T+\alpha \rbrack ,L^{2}(B))}\leq c(n,\epsilon )R^{\frac{n}{2}-\frac{n}{r}}{\left\Vert{\omega
 }\right\Vert}_{L^{s}(\lbrack 0,T+\alpha \rbrack ,L^{r}(B))}$\ \par 
with $c$ depending only on $n$ and $\epsilon .$ $\blacksquare $\ \par 
\ \par 
For $t\geq 0$ define $k:={\left\lceil{t}\right\rceil}\in {\mathbb{N}}$
 the integral part by excess, i.e.: $t\leq k<t+1.$\ \par 
Now set $k:={\left\lceil{\frac{n(r-2)}{2mr}}\right\rceil}$ then
 $k$ is the smallest integer such that, with $\frac{1}{r_{k}}=\frac{1}{2}-\frac{mk}{n},$
 we have $r_{k}\geq r.$\ \par 

\begin{prop}
~\label{pBB15}Let $r\geq 2.$ Let $B:=B(x,R)$ be a $\displaystyle
 (m,\epsilon )$-admissible ball and set $B^{k+1}:=B(x,2^{-k-1}R).$
 Then for any $\alpha >0$ we have the estimates, using $\beta
 ,\gamma ,\delta $ from Definition~\ref{pBB29}:\par 
\quad \quad \quad $\displaystyle R^{\delta }{\left\Vert{\partial _{t}u}\right\Vert}_{L^{s}(\lbrack
 0,T\rbrack ,L^{r_{k}}_{G}(B^{k+1}))}+R^{\gamma }{\left\Vert{u}\right\Vert}_{L^{s}(\lbrack
 0,T\rbrack ,W^{m,r_{k}}_{G}(B^{k+1}))}\leq $\par 
\quad \quad \quad \quad \quad $\displaystyle \leq c_{1}(k)R^{\beta }{\left\Vert{\omega }\right\Vert}_{L^{s}(\lbrack
 0,T+\alpha \rbrack ,L^{r}_{G}(B))}+c_{2}(k){\left\Vert{u}\right\Vert}_{L^{s}(\lbrack
 0,T+\alpha \rbrack ,L^{2}_{G}(B))},$\par 
and the constants $c_{1}(k),\ c_{2}(k)$ being independent of
 $B\in {\mathcal{A}}_{m}(\epsilon ).$\par 
In the case of functions instead of sections of $G$ we have the
 same estimates but with $B\in {\mathcal{A}}_{m-1}(\epsilon )$
 and using $\beta ',\gamma ',\delta '$ from Definition~\ref{pBB29}
 instead of $\displaystyle \beta ,\gamma ,\delta .$
\end{prop}
\quad Proof.\ \par 
Take $B:=B(x,R),\ B^{k}:=B(x,2^{-k}R).$\ \par 
\quad By the LIR inequality, Corollary~\ref{pBB11}, we get, with $\tau =2$:\ \par 
\quad \quad \quad $\displaystyle R^{m}{\left\Vert{\partial _{t}u}\right\Vert}_{L(2,1)}+R^{2m}{\left\Vert{u}\right\Vert}_{W(2,1)}\leq
 c_{3}(0)R^{m}{\left\Vert{Du}\right\Vert}_{L(2,0)}+c_{4}(0){\left\Vert{u}\right\Vert}_{L(2,0)},$\
 \par 
Now using Lemma~\ref{pBB14}, we get:\ \par 
\quad \quad \quad $\displaystyle {\left\Vert{\omega }\right\Vert}_{L(2,0)}\leq
 c(n,\epsilon )R^{\frac{n}{2}-\frac{n}{r}}{\left\Vert{\omega
 }\right\Vert}_{L(r,0)}.$\ \par 
Putting it above with $Du=\omega ,$ we get:\ \par 
\quad $\displaystyle R^{m}{\left\Vert{\partial _{t}u}\right\Vert}_{L(2,1)}+R^{2m}{\left\Vert{u}\right\Vert}_{W(2,1)}\leq
 c_{3}(0)c(n,\epsilon )R^{m}R^{\frac{n}{2}-\frac{n}{r}}{\left\Vert{\omega
 }\right\Vert}_{L(r,0)}+c_{4}(0){\left\Vert{u}\right\Vert}_{L(2,0)}.$\ \par 
Hence we have the induction hypothesis at level $k=0,$\ \par 
\quad $\displaystyle IH(0)\ \ \ \ \ \ R^{d_{0}}{\left\Vert{\partial
 _{t}u}\right\Vert}_{L(2,1)}+R^{b_{0}}{\left\Vert{u}\right\Vert}_{W(2,1)}\leq
 c_{1}(0)R^{a_{0}}{\left\Vert{\omega }\right\Vert}_{L(r,0)}+c_{2}(0){\left\Vert{u}\right\Vert}_{L(2,0)},$\
 \par 
with\ \par 
\quad \quad \quad $\displaystyle d_{0}=m,\ b_{0}=2m,\ a_{0}=m+\frac{n}{2}-\frac{n}{r}$\ \par 
and\ \par 
\quad \quad \quad $\displaystyle c_{1}(0)=c(n,\epsilon )c_{3}(0),\ c_{2}(0)=c_{4}(0).$\ \par 
\ \par 
So applying the induction Lemma~\ref{pL30}, we get\ \par 
\quad $\displaystyle IH(1)\ \ \ \ \ R^{d_{1}}{\left\Vert{\partial _{t}u}\right\Vert}_{L(\tau
 ,2)}+R^{b_{k+1}}{\left\Vert{u}\right\Vert}_{W(\tau ,2)}\leq
 c_{1}(1)R^{a_{1}}{\left\Vert{\omega }\right\Vert}_{L(r,0)}+c_{2}(1){\left\Vert{u}\right\Vert}_{L(2,0)}.$\
 \par 
with $\frac{1}{r_{1}}=\frac{1}{2}-\frac{m}{n},\ \tau :=\min (r_{1},r),$
 and\ \par 
\quad \quad \quad $\displaystyle d_{1}=3m+b_{0};\ b_{1}=4m+b_{0};\ a_{1}=\min (a_{0},3m+b_{0}),$\
 \par 
hence\ \par 
\quad \quad \quad $\displaystyle d_{1}=5m;\ b_{1}=6m;\ a_{1}=\min (m+\frac{n}{2}-\frac{n}{r},5m),$\
 \par 
and\ \par 
\quad \quad \quad $\displaystyle c_{1}(k+1)=c_{3}(k)+cc_{4}(k)c_{1}(k);\ c_{2}(k+1)=cc_{4}(k)c_{2}(k).$\
 \par 
By induction, we get\ \par 
\quad \quad \quad $\displaystyle b_{k}=4mk+2m;\ d_{k}=4mk+m,$\ \par 
and\ \par 
\quad \quad \quad $\displaystyle a_{0}=m+\frac{n}{2}-\frac{n}{r},\ \ \forall k\geq
 1,\ a_{k}=\min (m+\frac{n}{2}-\frac{n}{r},\ 5m).$\ \par 
\ \par 
\quad $\bullet $ if $r_{1}\geq r\Rightarrow \tau =r$  and we get:\ \par 
\quad \quad \quad $\displaystyle R^{d_{1}}{\left\Vert{\partial _{t}u}\right\Vert}_{L(r,2)}+R^{b_{1}}{\left\Vert{u}\right\Vert}_{W(r,2)}\leq
 c_{1}(1)R^{a_{1}}{\left\Vert{\omega }\right\Vert}_{L(r,0)}+c_{2}(1){\left\Vert{u}\right\Vert}_{L(2,0)}.$\
 \par 
And we are done.\ \par 
\ \par 
\quad $\bullet $ if $\tau =r_{1}<r,$ by the induction Lemma~\ref{pL30},
 after $k$ steps, we get with  $\frac{1}{r_{k}}=\frac{1}{2}-\frac{mk}{n},\
 \tau :=\min (r_{k},r)$:\ \par 
\quad \quad \quad $\displaystyle a_{k}=\min (m+\frac{n}{2}-\frac{n}{r},\ 5m),\
 b_{k}=(4k+2)m;\ d_{k}=(4k+1)m.$\ \par 
Then\ \par 
$\displaystyle IH(k)\ \ \ \ \ R^{d_{k}}{\left\Vert{\partial _{t}u}\right\Vert}_{L(\tau
 ,k+1)}+R^{b_{k}}{\left\Vert{u}\right\Vert}_{W(\tau ,k+1)}\leq
 c_{1}(k)R^{a_{k}}{\left\Vert{\omega }\right\Vert}_{L(r,0)}+c_{2}(k){\left\Vert{u}\right\Vert}_{L(2,0)}.$\
 \par 
\quad Hence if $r_{k}\geq r$ we are done as above, if not we repeat
 the process. Because $\frac{1}{r_{k}}=\frac{1}{2}-\frac{mk}{n}$
 after a finite number $k={\left\lceil{\frac{n(r-2)}{2mr})}\right\rceil}$
 of steps we have $r_{k}\geq r$ and we get, with $B^{k}:=B(x,R/2^{k})$:\ \par 
\quad $\displaystyle R^{d_{k}}{\left\Vert{\partial _{t}u}\right\Vert}_{L(r_{k},k+1)}+R^{b_{k}}{\left\Vert{u}\right\Vert}_{W(r_{k},k+1)}\leq
 c_{1}(k)R^{a_{k}}{\left\Vert{\omega }\right\Vert}_{L(r,0)}+c_{2}(k){\left\Vert{u}\right\Vert}_{L(2,0)}.$\
 \par 
Replacing the values of $L(r,k)$ and $W(r,k)$:\ \par 
\quad $\displaystyle R^{d_{k}}{\left\Vert{\partial _{t}u}\right\Vert}_{L^{s}(\lbrack
 0,T+\alpha /2^{k+1}\rbrack ,L^{r_{k}}_{G}(B^{k+1}))}+R^{b_{k}}{\left\Vert{u}\right\Vert}_{L^{s}(\lbrack
 0,T+\alpha /2^{k+1}\rbrack ,W^{m,r_{k}}_{G}(B^{k+1}))}\leq $\ \par 
\quad \quad \quad $\displaystyle \leq c_{1}(k)R^{a_{k}}{\left\Vert{\omega }\right\Vert}_{L^{s}(\lbrack
 0,T+\alpha \rbrack ,L^{r}_{G}(B))}+c_{2}(k){\left\Vert{u}\right\Vert}_{L^{s}(\lbrack
 0,T+\alpha \rbrack ,L^{2}_{G}(B))}.$\ \par 
With $c_{j}(k)$ depending on $\epsilon ,n,m,\alpha ,k$ and not on $B.$\ \par 
Because:\ \par 
\quad \quad \quad $\displaystyle {\left\Vert{\partial _{t}u}\right\Vert}_{L^{s}(\lbrack
 0,T\rbrack ,L^{r_{k}}_{G}(B^{k+1}))}\leq {\left\Vert{\partial
 _{t}u}\right\Vert}_{L^{s}(\lbrack 0,T+\alpha /2^{k+1}\rbrack
 ,L^{r_{k}}_{G}(B^{k+1}))}$\ \par 
and\ \par 
\quad \quad \quad $\displaystyle {\left\Vert{u}\right\Vert}_{L^{s}(\lbrack 0,T\rbrack
 ,W^{m,r_{k}}_{G}(B^{k+1}))}\leq {\left\Vert{u}\right\Vert}_{L^{s}(\lbrack
 0,T+\alpha /2^{k+1}\rbrack ,W^{m,r_{k}}_{G}(B^{k+1}))},$\ \par 
this proves the proposition for sections of $G.$\ \par 
\ \par 
\quad In the case of functions instead of sections of $G$ we have the
 same estimates but with $B\in {\mathcal{A}}_{m-1}(\epsilon )$ and:\ \par 
\quad $\displaystyle \forall k\geq 1,\ a_{k}=\min (m+\frac{n}{2}-\frac{n}{r},\
 4m-1),\ b_{k}=k(4m-1)+2m;\ d_{k}=m+k(4m-1),$\ \par 
and the constants $c_{1}(k),\ c_{2}(k)$ being independent of
 $B\in {\mathcal{A}}_{m-1}(\epsilon ).$\ \par 
This justifies the notation in Definition~\ref{pBB29}.\ \par 
\quad The proof is complete. $\blacksquare $\ \par 

\section{Vitali covering.}

\begin{lem}
Let ${\mathcal{F}}$ be a collection of balls $\lbrace B(x,r(x))\rbrace
 $ in a metric space, with $\forall B(x,r(x))\in {\mathcal{F}},\
 0<r(x)\leq R.$ There exists a disjoint subcollection ${\mathcal{G}}$
 of ${\mathcal{F}}$ with the following properties:\par 
\quad \quad every ball $B$ in ${\mathcal{F}}$ intersects a ball $C$ in ${\mathcal{G}}$
 and $B\subset 5C.$
\end{lem}
This is a well known lemma, see for instance ~\cite{EvGar92},
 section 1.5.1.\ \par 
\ \par 
\quad Fix $\epsilon >0$ and let $\forall x\in M,\ r(x):=R_{\epsilon
 }(x)/5,\ $where $R_{\epsilon }(x)$ is the admissible  radius
 at $\displaystyle x,$ we built a Vitali covering with the collection
 ${\mathcal{F}}:=\lbrace B(x,r(x))\rbrace _{x\in M}.$ The previous
 lemma gives a disjoint subcollection ${\mathcal{G}}$ such that
 every ball $B$ in ${\mathcal{F}}$ intersects a ball $C$ in ${\mathcal{G}}$
 and we have $\displaystyle B\subset 5C.$ We set ${\mathcal{D}}(\epsilon
 ):=\lbrace x\in M::B(x,r(x))\in {\mathcal{G}}\rbrace $ and ${\mathcal{C}}_{\epsilon
 }:=\lbrace B(x,5r(x)),\ x\in {\mathcal{D}}(\epsilon )\rbrace
 $: we shall call ${\mathcal{C}}_{\epsilon }$ a $\epsilon $-\textbf{admissible
 covering} of $(M,g).$\ \par 
\quad More generally  let $k\in {\mathbb{N}}$ and consider the collection
 $\displaystyle {\mathcal{F}}_{k}(\epsilon ):=\lbrace B(x,r_{k}(x))\rbrace
 _{x\in M}$ where, for $x\in M,\ r_{k}(x):=2^{-k}R_{\epsilon
 }(x)/5\eta ,\ $ still where $R_{\epsilon }(x)$ is the admissible
  $\epsilon $-radius at $\displaystyle x.$ The integer $\eta
 \geq 1$ will be chosen later. The previous lemma gives a disjoint
 subcollection ${\mathcal{G}}_{k}(\epsilon )$ such that every
 ball $B$ in ${\mathcal{F}}_{k}(\epsilon )$ intersects a ball
 $C$ in ${\mathcal{G}}_{k}(\epsilon )$ and we have $\displaystyle
 B\subset 5C.$ We set ${\mathcal{D}}_{k}(\epsilon ):=\lbrace
 x\in M::B(x,r_{k}(x))\in {\mathcal{G}}_{k}(\epsilon )\rbrace
 $ and ${\mathcal{C}}_{k}(\epsilon ):=\lbrace B(x,5r_{k}(x)),\
 x\in {\mathcal{D}}_{k}(\epsilon )\rbrace $: we shall call ${\mathcal{C}}_{k}(\epsilon
 )$ a $(k,\epsilon )$-\textbf{admissible covering} of $(M,g).$\ \par 
\quad We have the lemma:\ \par 

\begin{lem}
~\label{pL47}Let $\displaystyle B(x,5r_{k}(x))\in {\mathcal{C}}_{k}(\epsilon
 )$ then $B^{0}(x,\tilde R(x))$ with $\displaystyle \tilde R(x):=2^{k}{\times}5r_{k}(x))$
 is still a $\epsilon $-admissible ball. Moreover we have that
 all the balls $B^{j}(x):=B^{0}(x,2^{-j}\tilde R(x)),\ j=0,1,...,\
 k$ are also $\displaystyle \epsilon $-admissible balls and $\lbrace
 B^{j}(x),\ x\in {\mathcal{D}}_{k}(\epsilon )\rbrace ,$ for $j=0,...,k,$
 is a covering of $M.$ 
\end{lem}
\quad Proof.\ \par 
Take $x\in {\mathcal{D}}_{k}(\epsilon )$ then we have that the
 geodesic ball $\displaystyle B^{0}(x,\tilde R(x))=B(x,R_{\epsilon
 }/\eta )$ is $\epsilon $-admissible and because $\displaystyle
 2^{-k}R_{\epsilon }/\eta <R_{\epsilon },$ for $\eta \geq 1,$
 we get that $\displaystyle B(x,R(x))$ is also $\epsilon $-admissible.\ \par 
\quad The same for $\displaystyle B^{j}(x)=B^{0}(x,2^{-j}\tilde R(x))$
 because $2^{-j}\tilde R(x)<2^{-j}R_{\epsilon }(x)/\eta .$\ \par 
The fact that $\lbrace B^{k}(x),\ x\in {\mathcal{D}}_{k}(\epsilon
 )\rbrace $ is a covering of $M$ is just the Vitali lemma and,
 because $j\leq k\Rightarrow B^{j}(x)\supset B^{k}(x),$ we get
 that $\lbrace B^{j}(x),\ x\in {\mathcal{D}}_{k}(\epsilon )\rbrace
 $ is also a covering of $M.$ $\blacksquare $\ \par 
\ \par 
\quad Then we have:\ \par 

\begin{prop}
~\label{CF2}Let $(M,g)$ be a riemannian manifold, then the overlap
 of a $(k,\epsilon )$-admissible covering ${\mathcal{C}}_{k}(\epsilon
 )$ is less than $T=\frac{(1+\epsilon )^{n/2}}{(1-\epsilon )^{n/2}}(100)^{n},$
 i.e.\par 
\quad \quad \quad $\forall x\in M,\ x\in B(y,5r_{k}(y))$ where $B(y,r_{k}(y))\in
 {\mathcal{G}}_{k}(\epsilon )$ for at most $T$ such balls.\par 
Moreover we have\par 
\quad \quad \quad $\displaystyle \forall f\in L^{1}(M),\ \sum_{j\in {\mathbb{N}}}{\int_{B(x_{j},r_{k}(x_{j}))}{\left\vert{f(x)}\right\vert
 dv_{g}(x)}}\leq T{\left\Vert{f}\right\Vert}_{L^{1}(M)}.$
\end{prop}
\quad Proof.\ \par 
Let $B_{j}:=B(x_{j},r_{k}(x_{j}))\in {\mathcal{G}}_{k}(\epsilon
 )$ and suppose that $x\in \bigcap_{j=1}^{l}{B(x_{j},5r_{k}(x_{j}))}.$
 Then we have\ \par 
\quad \quad \quad $\displaystyle \forall j=1,...,l,\ d(x,x_{j})\leq 5r_{k}(x_{j}).$\ \par 
Hence\ \par 
\quad $\displaystyle d(x_{j},x_{m})\leq d(x_{j},x)+d(x,x_{m})\leq 5(r_{k}(x_{j})+r_{k}(x_{m}))\leq
 2^{-k}(R_{\epsilon }(x_{j})+R_{\epsilon }(x_{m}))/\eta .$\ \par 
Suppose that $r_{k}(x_{j})\geq r_{k}(x_{m})$ then $x_{m}\in B(x_{j},10r_{k}(x_{j}))\subset
 B(x_{j},R_{\epsilon }(x_{j}))$ because $10r_{k}(x_{j})=10{\times}2^{-k}R_{\epsilon
 }(x_{j})/\eta \leq R_{\epsilon }(x_{j})$ because now on \emph{we
 choose } $\eta =10.$ Then, by the slow variation of the $\epsilon
 $-radius Lemma~\ref{m7}, we have $R_{\epsilon }(x_{j})\leq 2R_{\epsilon
 }(x_{m}).$ If $\displaystyle r_{k}(x_{j})\leq r_{k}(x_{m})$
 then, the same way, $\displaystyle 2R_{\epsilon }(x_{j})\geq
 R_{\epsilon }(x_{m}).$ Hence in any case we have $\displaystyle
 \frac{1}{2}R_{\epsilon }(x_{m})\leq R_{\epsilon }(x_{j})\leq
 2R_{\epsilon }(x_{m}).$\ \par 
So $\displaystyle d(x_{j},x_{m})\leq 2^{-k}{\times}3R_{\epsilon
 }(x_{j})/10$ hence\ \par 
\quad $\displaystyle \forall m=1,...,l,\ B(x_{m},\ r_{k}(x_{m}))\subset
 B(x_{j},\ 2^{-k}{\times}3R_{\epsilon }(x_{j})/10+2^{-k}R_{\epsilon
 }(x_{m})/10)\subset B(x_{j},\ 2^{-k}{\times}5R_{\epsilon }(x_{j})/10)$\ \par 
because $\displaystyle R_{\epsilon }(x_{m})\leq 2R_{\epsilon }(x_{j}).$\ \par 
The balls in ${\mathcal{G}}_{k}(\epsilon )$ being disjoint, we
 get, setting $B_{m}:=B(x_{m},\ r_{k}(x_{m})),$\ \par 
\quad \quad \quad $\displaystyle \sum_{m=1}^{l}{\mathrm{V}\mathrm{o}\mathrm{l}(B_{m})}\leq
 \mathrm{V}\mathrm{o}\mathrm{l}(B(x_{j},2^{-k}{\times}5R_{\epsilon
 }(x_{j})/\eta )=\mathrm{V}\mathrm{o}\mathrm{l}(B(x_{j},5r_{k}(x_{j})).$\ \par 
\quad The Lebesgue measure read in the chart $\varphi $  and the canonical
 measure $dv_{g}$ on $B(x,R_{\epsilon }(x))$ are equivalent;
 precisely because of condition 1) in the admissible ball definition,
 we get that:\ \par 
\quad \quad \quad $\displaystyle (1-\epsilon )^{n}\leq \left\vert{\mathrm{d}\mathrm{e}\mathrm{t}g}\right\vert
 \leq (1+\epsilon )^{n},$\ \par 
and the measure $dv_{g}$ read in the chart $\varphi $ is $dv_{g}={\sqrt{\left\vert{\mathrm{d}\mathrm{e}\mathrm{t}g_{ij}}\right\vert
 }}d\xi ,$ where $\displaystyle d\xi $ is the Lebesgue measure
 in ${\mathbb{R}}^{n}.$ In particular:\ \par 
\quad \quad \quad $\displaystyle \forall x\in M,\ \mathrm{V}\mathrm{o}\mathrm{l}(B(x,\
 R_{\epsilon }(x)))\leq (1+\epsilon )^{n/2}\nu _{n}R^{n},$\ \par 
where $\nu _{n}$ is the euclidean volume of the unit ball in
 ${\mathbb{R}}^{n}.$\ \par 
\quad Now because $R_{\epsilon }(x_{j})$ is the admissible radius and
 $5r_{k}(x_{j})=2^{-k}{\times}5R_{\epsilon }(x_{j})/10<R_{\epsilon
 }(x_{j}),$ because $\eta =10,$\ \par 
\quad \quad \quad $\displaystyle \mathrm{V}\mathrm{o}\mathrm{l}(B(x_{j},5r_{k}(x_{j})))\leq
 5^{n}(1+\epsilon )^{n/2}v_{n}r_{k}(x_{j})^{n}.$\ \par 
On the other hand we have also\ \par 
\quad \quad \quad $\displaystyle \mathrm{V}\mathrm{o}\mathrm{l}(B_{m})\geq v_{n}(1-\epsilon
 )^{n/2}r_{k}(x_{m})^{n}\geq v_{n}(1-\epsilon )^{n/2}2^{-n}r_{k}(x_{j})^{n},$\
 \par 
hence\ \par 
\quad \quad \quad $\displaystyle \sum_{j=1}^{l}{(1-\epsilon )^{n/2}2^{-n}r(x_{j})^{n}}\leq
 5^{n}(1+\epsilon )^{n/2}r_{k}(x_{j})^{n},$\ \par 
so finally\ \par 
\quad \quad \quad $\displaystyle l\leq (5{\times}2)^{n}\frac{(1+\epsilon )^{n/2}}{(1-\epsilon
 )^{n/2}},$\ \par 
which means that $T\leq \frac{(1+\epsilon )^{n/2}}{(1-\epsilon
 )^{n/2}}(100)^{n}.$\ \par 
\quad Saying that any $x\in M$ belongs to at most $T$ balls of the
 covering $\lbrace B_{j}\rbrace $ means that $\ \sum_{j\in {\mathbb{N}}}{{\11}_{B_{j}}(x)}\leq
 T,$ and this implies easily that:\ \par 
\quad \quad \quad $\displaystyle \forall f\in L^{1}(M),\ \sum_{j\in {\mathbb{N}}}{\int_{B_{j}}{\left\vert{f(x)}\right\vert
 dv_{g}(x)}}\leq T{\left\Vert{f}\right\Vert}_{L^{1}(M)}.$\ \par 
The proof is complete. $\blacksquare $\ \par 

\begin{cor}
~\label{pL48} Let $(M,g)$ be a riemannian manifold. Consider
 the $(k,\epsilon )$-admissible covering ${\mathcal{C}}_{k}(\epsilon
 ).$ Then the overlap of the associated covering by the balls
 $\lbrace B^{0}(x,R_{\epsilon }(x)/\eta )\rbrace _{x\in {\mathcal{D}}_{k}(\epsilon
 )}$ verifies $\displaystyle T_{k}\leq \frac{(1+\epsilon )^{n/2}}{(1-\epsilon
 )^{n/2}}(100)^{n}{\times}2^{nk}.$
\end{cor}
\quad Proof.\ \par 
We start the proof exactly the same way as above.\ \par 
Let $B_{j}:=B(x_{j},R_{\epsilon }(x_{j})/\eta ),\ x_{j}\in {\mathcal{D}}_{k}(\epsilon
 )$ and suppose that $x\in \bigcap_{j=1}^{l}{B(x_{j},R_{\epsilon
 }(x_{j})/\eta )}.$ Then we have, as above:\ \par 
\quad \quad \quad $\displaystyle \forall m=1,...,l,\ B(x_{m},\ R_{\epsilon }(x_{m}))\subset
 B(x_{j},\ 3R_{\epsilon }(x_{j})/\eta +R_{\epsilon }(x_{m})/\eta
 )\subset B(x_{j},\ 5R_{\epsilon }(x_{j})/\eta ).$\ \par 
The balls in ${\mathcal{G}}_{k}(\epsilon )$ being disjoint, we
 get, setting $B_{m}:=B(x_{m},\ r_{k}(x_{m})),$\ \par 
\quad $\displaystyle \sum_{m=1}^{l}{\mathrm{V}\mathrm{o}\mathrm{l}(B_{m})}\leq
 \mathrm{V}\mathrm{o}\mathrm{l}(B(x_{j},5R_{\epsilon }(x_{j})/\eta
 )=\mathrm{V}\mathrm{o}\mathrm{l}(B(x_{j},5{\times}2^{k}r_{k}(x_{j})).$\ \par 
Exactly as above, we get because of the factor $2^{k},$\ \par 
\quad \quad \quad $\displaystyle l\leq (5{\times}2{\times}2^{k})^{n}\frac{(1+\epsilon
 )^{n/2}}{(1-\epsilon )^{n/2}},$\ \par 
hence the result. $\blacksquare $\ \par 

\section{The threshold.~\label{p32}}
\quad We shall need the following "threshold hypothesis".\ \par 
(THL2)  For any $\omega \in L^{s}(\lbrack 0,T\rbrack ,L^{2}_{G}(M))$
 there is a $\displaystyle u\in L^{s}(\lbrack 0,T\rbrack ,L^{2}_{G}(M))$
 such that $Du=\omega $ with the estimate:\ \par 
\quad \quad \quad $\displaystyle {\left\Vert{u}\right\Vert}_{L^{s}(\lbrack 0,T\rbrack
 ,L^{2}_{G}(M))}\lesssim {\left\Vert{\omega }\right\Vert}_{L^{s}(\lbrack
 0,T\rbrack ,L^{2}_{G}(M))}.$\ \par 
\quad We have that the (THL2) hypothesis is true for the heat equation.\ \par 
My first proof used Patodi Hodge decomposition~\cite{Patodi71}
 in compact manifold, but E-M. Ouhabaz gave me the following
 one, see~\cite{MagOuh17}, page 2, which is shorter.\ \par 
\quad This is the only place where we use estimates on the heat semi-group
 $\displaystyle (e^{-t\Delta })_{t\geq 0}$ on the manifold. Here
 we set $G:=\Lambda ^{p}(M),$ the bundle of $p$-forms on $M.$\ \par 

\begin{thm}
~\label{pBB19}Let $M$ be a connected complete ${\mathcal{C}}^{2}$
 riemannian manifold. Let, for $\displaystyle t\geq 0,\ \omega
 (t,x)\in L^{2}_{p}(M).$ Then we have a solution $u$ of the heat
 equation $\displaystyle Du:=\partial _{t}u+\Delta u=\omega ,\
 u(0,x)\equiv 0,$ such that $\displaystyle \forall t\geq 0,\
 u(t,x)\in L^{2}_{p}(M)$ with the estimate:\par 
\quad \quad \quad $\displaystyle \forall t\geq 0,\ {\left\Vert{u(t,\cdot )}\right\Vert}_{L^{2}_{p}(M)}\leq
 \int_{0}^{t}{{\left\Vert{\omega (\tau ,\cdot )}\right\Vert}_{L^{2}_{p}(M)}d\tau
 }.$
\end{thm}
\quad Proof.\ \par 
It is well known that the Hodge laplacian is essentially positive
 on $p$-forms in $L^{2}_{p}(M),$ so $(e^{-t\Delta })_{t\geq 0}$
 is a contraction semi-group on $\displaystyle L^{2}_{p}(M).$
 Take $p(t,x,y)$ the kernel associated to the semi-group $\displaystyle
 (e^{-t\Delta })_{t\geq 0}.$ We have the non homogeneous solution:\ \par 
\quad \quad \quad $\displaystyle u(t,x)=\int_{0}^{t}{\int_{M}{p(t-s,x,y)\omega
 (s,y)dv(y)ds}}$\ \par 
which verifies $\displaystyle Du:=\partial _{t}u+\Delta u=\omega
 ,\ u(0,x)\equiv 0.$\ \par 
Fix $s,t$ and set:\ \par 
\quad \quad \quad $\displaystyle v(t,s,x):=\int_{M}{p(t-s,x,y)\omega (s,y)dv(y)}$\ \par 
then the contraction property gives:  ${\left\Vert{v(t,s,\cdot
 )}\right\Vert}_{L^{2}_{p}(M)}\leq {\left\Vert{\omega (s,\cdot
 )}\right\Vert}_{L^{2}_{p}(M)}.$  Hence: ${\left\Vert{u(t,\cdot
 )}\right\Vert}_{L^{2}_{p}(M)}\leq \int_{0}^{t}{{\left\Vert{\omega
 (s,\cdot )}\right\Vert}_{L^{2}_{p}(M)}ds.}$\ \par 
The proof is complete. $\blacksquare $\ \par 
\quad Clearly this result implies the hypothesis (THL2).\ \par 

\section{Global results.}
\quad We want to globalise Theorem~\ref{m6} by use of our Vitali coverings. \ \par 

\begin{lem}
~\label{m8}We have for any section $f:\ M\rightarrow G$ and $\tau
 \in (1,\infty ),$ with $w(x):=R_{\epsilon }(x)^{\gamma \tau
 }$ and $B(x):=B(x,R_{x}(x)/10),\ B^{k}(x):=B(x,2^{-k}R_{\epsilon
 }(x)/10),$ that:\par 
\quad \quad \quad $\displaystyle \forall \tau \geq 1,\ {\left\Vert{f}\right\Vert}_{W^{l,\tau
 }_{G}(M,\ w)}^{\tau }\simeq \sum_{x\in {\mathcal{D}}_{k}(\epsilon
 )}{R_{\epsilon }(x)^{\gamma \tau }{\left\Vert{f}\right\Vert}_{W^{l,\tau
 }_{G}(B^{k}(x))}^{\tau }};$\par 
and:\par 
\quad \quad \quad $\displaystyle \forall \tau \geq 1,\ {\left\Vert{f}\right\Vert}_{W^{l,\tau
 }_{G}(M,w)}^{\tau }\simeq \sum_{x\in {\mathcal{D}}_{k}(\epsilon
 )}{R_{\epsilon }(x)^{\gamma \tau }{\left\Vert{f}\right\Vert}_{W^{l,\tau
 }_{G}(B(x))}^{\tau }}.$
\end{lem}
\quad Proof.\ \par 
Let $x\in {\mathcal{D}}_{k}(\epsilon ),$ this implies that $B^{k}(x):=B(x,2^{-k}R_{\epsilon
 }(x)/10)\in {\mathcal{C}}_{k}(\epsilon ).$ \ \par 
\quad $\bullet $ First we start with $l=0.$ We shall deal with the
 function $\left\vert{f}\right\vert .$\ \par 
\quad We have, because ${\mathcal{C}}_{k}(\epsilon )$ is a covering
 of $M$ and with $\forall y\in B(x),\ R(y)\ :=R_{\epsilon }\ (y)$\ \par 
\quad \quad \quad $\displaystyle {\left\Vert{f}\right\Vert}_{L^{\tau }(M,w)}^{\tau
 }:=\int_{M}{\left\vert{f(x)}\right\vert ^{\tau }w(x)dv(x)}\leq
 \sum_{x\in {\mathcal{D}}_{k}(\epsilon )}{\int_{B^{k}(x)}{\left\vert{f(y)}\right\vert
 ^{\tau }R(y)^{\gamma \tau }}dv(y)}.$\ \par 
We have, by Lemma~\ref{m7}, $\forall y\in B,\ R(y)\leq 2R(x),$ then\ \par 
\quad \quad \quad $\displaystyle \sum_{x\in {\mathcal{D}}_{k}(\epsilon )}{\int_{B^{k}(x)}{\left\vert{f(y)}\right\vert
 ^{\tau }R(y)^{\gamma \tau }}dv(y)}\leq $\ \par 
\quad \quad \quad \quad \quad $\displaystyle \leq \sum_{x\in {\mathcal{D}}_{k}(\epsilon )}{2^{\gamma
 \tau }R(x)^{\gamma \tau }\int_{B^{k}(x)}{\left\vert{f(y)}\right\vert
 ^{\tau }}dv(y)}\leq 2^{\gamma \tau }\sum_{x\in {\mathcal{D}}_{k}(\epsilon
 )}{R(x)^{\gamma \tau }{\left\Vert{f}\right\Vert}_{L^{\tau }(B^{k}(x))}^{\tau
 }}.$\ \par 
Hence\ \par 
\quad \quad \quad $\displaystyle {\left\Vert{f}\right\Vert}_{L_{G}^{\tau }(M,w)}^{\tau
 }\leq 2^{\gamma \tau }\sum_{x\in {\mathcal{D}}_{k}(\epsilon
 )}{R(x)^{\gamma \tau }{\left\Vert{f}\right\Vert}_{L_{G}^{\tau
 }(B^{k})}^{\tau }}.$\ \par 
\ \par 
\quad To get the converse inequality we still use Lemma~\ref{m7}: $\forall
 y\in B,\ R(x)\leq 2R(y)$ so we get:\ \par 
\quad $\displaystyle \sum_{x\in {\mathcal{D}}_{k}(\epsilon )}{R(x)^{\gamma
 \tau }\int_{B^{k}(x)}{\left\vert{f(y)}\right\vert ^{\tau }}dv(y)}\leq
 2^{\gamma \tau }\sum_{x\in {\mathcal{D}}_{k}(\epsilon )}{\int_{B^{k}(x)}{R(y)^{\gamma
 \tau }\left\vert{f(y)}\right\vert ^{\tau }}dv(y)}.$\ \par 
Now we use the fact that the overlap of ${\mathcal{C}}_{k}(\epsilon
 )$ is bounded by $T,$\ \par 
\quad $\displaystyle \sum_{x\in {\mathcal{D}}_{k}(\epsilon )}{\int_{B^{k}(x)}{R(y)^{\gamma
 \tau }\left\vert{f(y)}\right\vert ^{\tau }}dv(y)}\leq 2^{\gamma
 \tau }T\int_{M}{R(y)^{\gamma \tau }\left\vert{f(y)}\right\vert
 ^{\tau }}dv(y)=2^{\gamma \tau }T{\left\Vert{f}\right\Vert}_{L^{\tau
 }(M,w)}^{\tau }.$\ \par 
So\ \par 
\quad \quad \quad $\displaystyle \sum_{x\in {\mathcal{D}}_{k}(\epsilon )}{R^{\gamma
 \tau }{\left\Vert{f}\right\Vert}_{L^{\tau }(B^{k})}}^{\tau }\leq
 2^{\gamma \tau }T{\left\Vert{f}\right\Vert}_{L^{\tau }(M,w)}^{\tau }.$\ \par 
\ \par 
\quad We already know that $\displaystyle \lbrace B::B^{k}\in {\mathcal{C}}_{k}(\epsilon
 )\rbrace $ is a covering of $M$ with a bounded overlap by Corollary~\ref{pL48},
 so we follow exactly the same lines to prove:\ \par 
\quad \quad \quad $\displaystyle \forall \tau \geq 1,\ {\left\Vert{f}\right\Vert}_{L^{\tau
 }_{G}(M,w)}^{\tau }\simeq \sum_{x\in {\mathcal{D}}_{k}(\epsilon
 )}{R(x)^{\gamma \tau }{\left\Vert{f}\right\Vert}_{L^{\tau }_{G}(B(x))}^{\tau
 }}.$\ \par 
\quad $\bullet $ Now let $l\geq 1.$\ \par 
We apply the case $l=0$ to the covariant derivatives of $f.$\ \par 
\quad \quad \quad $\displaystyle \forall \tau \geq 1,\ {\left\Vert{\nabla ^{l}f}\right\Vert}_{L^{\tau
 }_{G}(M,w)}^{\tau }\simeq \sum_{x\in {\mathcal{D}}_{k}(\epsilon
 )}{R(x)^{\gamma \tau }{\left\Vert{\nabla ^{l}f}\right\Vert}_{L^{\tau
 }_{G}(B(x))}^{\tau }}.$\ \par 
Because\ \par 
\quad \quad \quad $\displaystyle {\left\Vert{f}\right\Vert}_{W^{l,\tau }}={\left\Vert{f}\right\Vert}_{L^{\tau
 }}+\cdot \cdot \cdot +{\left\Vert{\nabla ^{l}f}\right\Vert}_{L^{\tau }}$\ \par 
we get\ \par 
\quad \quad \quad $\displaystyle \forall \tau \geq 1,\ {\left\Vert{\nabla ^{l}f}\right\Vert}_{W^{l,\tau
 }_{G}(M,w)}^{\tau }\simeq \sum_{x\in {\mathcal{D}}_{k}(\epsilon
 )}{R(x)^{\gamma \tau }{\left\Vert{\nabla ^{l}f}\right\Vert}_{W^{l,\tau
 }_{G}(B(x))}^{\tau }}.$\ \par 
The proof is complete. $\blacksquare $\ \par 

\begin{thm}
~\label{pBB16}Suppose that $A$ is a $(C,\theta )$-elliptic operator
 of order $m$ acting on sections of the adapted vector bundle
 $\displaystyle G:=(H,\pi ,M)$ in the complete riemannian manifold
 $(M,g),$ with $\theta <\pi /2,$ and consider the parabolic equation
 $\displaystyle Du=\partial _{t}u-Au$ also acting on sections of $G.$\par 
Let $\displaystyle R(x)=R_{m,\epsilon }(x)$ be the $(m,\epsilon
 )$ radius at the point $x\in M.$ Set $\displaystyle w_{1}(x):=R(x)^{\delta
 },\ w_{2}(x):=R(x)^{\gamma },\ w_{3}(x):=R(x)^{\beta },$ with
 the notation in Definition~\ref{pBB29}. We have:\par 
\quad \quad \quad $\displaystyle {\left\Vert{\partial _{t}u}\right\Vert}_{L^{r}(\lbrack
 0,T\rbrack ,L_{G}^{r}(M,w_{1}))}+{\left\Vert{u}\right\Vert}_{L^{r}(\lbrack
 0,T\rbrack ,W_{G}^{m,r}(M,w_{2}))}\leq $\par 
\quad \quad \quad \quad \quad $\displaystyle \leq c_{1}{\left\Vert{Du}\right\Vert}_{L^{r}(\lbrack
 0,T+\alpha \rbrack ,L_{G}^{r}(M,w_{3}))}+c_{2}{\left\Vert{u}\right\Vert}_{L^{r}(\lbrack
 0,T+\alpha \rbrack ,L_{G}^{2}(M))}.$\par 
In the case of functions instead of sections of $G$ we have the
 same estimates but with $R(x)=R_{m-1,\epsilon }(x)$ and:\par 
\quad \quad \quad $\displaystyle w_{1}(x):=R(x)^{\delta '},\ w_{2}(x):=R(x)^{\gamma
 '},\ w_{3}(x):=R(x)^{\beta '}.$
\end{thm}
\quad Proof.\ \par 
Once again we shall use the notation, for $k\geq 1,$\ \par 
\quad \quad \quad $\displaystyle L(s,r,k):=L^{s}(\lbrack 0,T\rbrack ,L^{r}_{G}(B^{k}));\
 \ W(s,r,k):=L^{s}(\lbrack 0,T\rbrack ,W^{m,r}_{G}(B^{k}))$\ \par 
and for $k=0,$\ \par 
\quad \quad \quad $\displaystyle L(s,r,0):=L^{s}(\lbrack 0,T+\alpha \rbrack ,L^{r}_{G}(B));\
 \ W(s,r,0):=L^{s}(\lbrack 0,T+\alpha \rbrack ,W^{m,r}_{G}(B))$\ \par 
With this notation, the Proposition~\ref{pBB15} gives, for $r\geq
 2$ and $\displaystyle k:={\left\lceil{\frac{n(r-2)}{2mr}}\right\rceil}$:\ \par 
\quad \quad \quad $\displaystyle R^{\delta }{\left\Vert{\partial _{t}u}\right\Vert}_{L(s,r_{k},k+1)}+R^{\gamma
 }{\left\Vert{u}\right\Vert}_{W(s,r_{k},k+1)}\leq c_{1}R^{\beta
 }{\left\Vert{\omega }\right\Vert}_{L(s,r,0)}+c_{2}{\left\Vert{u}\right\Vert}_{L(s,2,0)}.$\
 \par 
Because $\displaystyle r_{k}\geq r$ we get\ \par 
\quad \quad \quad $\displaystyle R^{\delta }{\left\Vert{\partial _{t}u}\right\Vert}_{L(s,r,k+1)}+R^{\gamma
 }{\left\Vert{u}\right\Vert}_{W(s,r,k+1)}\leq c_{1}R^{\beta }{\left\Vert{\omega
 }\right\Vert}_{L(s,r,0)}+c_{2}{\left\Vert{u}\right\Vert}_{L(s,2,0)}.$\ \par 
So for $s=r$ we get\ \par 
\quad \quad \quad \begin{equation} R^{\delta }{\left\Vert{\partial _{t}u}\right\Vert}_{L(r,r,k+1)}+R^{\gamma
 }{\left\Vert{u}\right\Vert}_{W(r,r,k+1)}\leq c_{1}R^{\beta }{\left\Vert{\omega
 }\right\Vert}_{L(r,r,0)}+c_{2}{\left\Vert{u}\right\Vert}_{L(r,2,0)}.\label{pBB30}\end{equation}\
 \par 
Because\ \par 
\quad \quad \quad $\displaystyle a+b\leq c+d\Rightarrow a^{r}+b^{r}\leq Ac^{r}+Bd^{r}$\ \par 
with constants $A,B$ depending on $r$ only, the inequality~(\ref{pBB30})
 can be read with a slight change of the constants:\ \par 
\quad \quad \quad $\displaystyle R^{r\delta }{\left\Vert{\partial _{t}u}\right\Vert}_{L(r,r,k+1)}^{r}+R^{r\gamma
 }{\left\Vert{u}\right\Vert}_{W(r,r,k+1)}^{r}\leq c_{1}(k)R^{r\beta
 }{\left\Vert{\omega }\right\Vert}_{L(r,r,0)}^{r}+c_{2}{\left\Vert{u}\right\Vert}_{L(r,2,0)}^{r}.$\
 \par 
By use of Lemma~\ref{m8} with $l=0,\ \tau =r,\ w_{1}(x):=R(x)^{r\delta
 },$\ \par 
\quad \quad \quad $\displaystyle \ {\left\Vert{\partial _{t}u}\right\Vert}_{L^{r}_{G}(M,\
 w_{1})}^{r}\simeq \sum_{x\in {\mathcal{D}}_{k+1}(\epsilon )}{R(x)^{\delta
 r}{\left\Vert{\partial _{t}u}\right\Vert}_{L^{r}_{G}(B^{k+1}(x))}^{r}};$\ \par 
hence integrating in $t\in \lbrack 0,T\rbrack $ we get:\ \par 
\quad \quad \quad $\displaystyle {\left\Vert{\partial _{t}u}\right\Vert}_{L^{r}(\lbrack
 0,T\rbrack ,L^{r}_{G}(M,w_{1}))}^{r}\simeq \sum_{x\in {\mathcal{D}}_{k+1}(\epsilon
 )}{R(x)^{\delta r}{\left\Vert{\partial _{t}u}\right\Vert}_{L^{r}(\lbrack
 0,T\rbrack ,L^{r}_{G}(B^{k+1}(x))}^{r}}.$\ \par 
The same way, with $\displaystyle l=m,\ \tau =r,\ w_{2}(x):=R(x)^{r\gamma
 },$\ \par 
\quad $\displaystyle {\left\Vert{u}\right\Vert}_{L^{r}(\lbrack 0,T\rbrack
 ,W^{m,r}_{G}(M,w_{2}))}^{r}\simeq \sum_{x\in {\mathcal{D}}_{k+1}}{R(x)^{\gamma
 r}{\left\Vert{u}\right\Vert}_{L^{r}(\lbrack 0,T\rbrack ,W^{m,r}_{G}(B^{k+1}(x)))}^{r}}$\
 \par 
with $l=0,\ \tau =r,\ w_{3}(x):=R(x)^{r\beta },$\ \par 
\quad $\displaystyle {\left\Vert{\omega }\right\Vert}_{L^{r}(\lbrack
 0,T+\alpha \rbrack ,L^{r}_{G}(M,w_{3}))}^{r}\simeq \sum_{x\in
 {\mathcal{D}}_{k+1}}{R(x)^{\beta r}{\left\Vert{\omega }\right\Vert}_{L^{r}(\lbrack
 0,T+\alpha \rbrack ,L^{r}_{G}(B(x)))}^{r}}$\ \par 
and with $l=0,\ \tau =r,$\ \par 
\quad \quad \quad $\displaystyle {\left\Vert{u}\right\Vert}_{L^{r}(\lbrack 0,T+\alpha
 \rbrack ,L^{2}_{G}(M))}^{r}\simeq \sum_{x\in {\mathcal{D}}_{k+1}}{{\left\Vert{u}\right\Vert}_{L^{r}(\lbrack
 0,T+\alpha \rbrack ,L^{2}_{G}(B(x)))}^{r}}.$\ \par 
\quad So, putting this in Theorem~\ref{m6}, we get, with $\displaystyle
 w_{1}(x):=R(x)^{r\delta },\ w_{2}(x):=R(x)^{r\gamma },\ w_{3}(x):=R(x)^{r\beta
 },$\ \par 
\quad \quad \quad $\displaystyle {\left\Vert{\partial _{t}u}\right\Vert}_{L^{r}(\lbrack
 0,T\rbrack ,L_{G}^{r}(M,w_{1}))}+{\left\Vert{u}\right\Vert}_{L^{r}(\lbrack
 0,T\rbrack ,W_{G}^{m,r}(M,w_{2}))}\leq $\ \par 
\quad \quad \quad \quad \quad $\displaystyle \leq c_{1}(k){\left\Vert{Du}\right\Vert}_{L^{r}(\lbrack
 0,T+\alpha \rbrack ,L_{G}^{r}(M,w_{3}))}+c_{2}(k){\left\Vert{u}\right\Vert}_{L^{r}(\lbrack
 0,T+\alpha \rbrack ,L_{G}^{2}(M))}.$\ \par 
The results, for functions instead of sections of $G,$ follow
 the same lines and we have the same estimates but with $R(x)=R_{m-1,\epsilon
 }(x)$ and the weights:\ \par 
\quad \quad \quad $\displaystyle w_{1}(x):=R(x)^{r\delta '},\ w_{2}(x):=R(x)^{r\gamma
 '},\ w_{3}(x):=R(x)^{r\beta '}.$\ \par 
The proof is complete. $\blacksquare $\ \par 

\begin{rem}
The weights $w_{j}(x)$ depend on $r$ but also on $m$ and $n$
 via $\beta ,\gamma $ and $\delta $ given by the Definition~\ref{pBB29}.
\end{rem}

\begin{cor}
~\label{pBB20}Let $M$ be a connected complete $n$-dimensional
 ${\mathcal{C}}^{m}$ riemannian manifold without boundary. Let
 $G:=(H,\pi ,M)$ be a complex ${\mathcal{C}}^{m}$ adapted vector
 bundle over $M.$ Suppose $\displaystyle Du:=\partial _{t}u-Au,$
 where $A$ is $(C,\theta )$-elliptic of order $m$ acting on sections
 of $G$ with $\theta <\pi /2.$ Moreover suppose we have (THL2).
 Let $r\geq 2$ and:\par 
\quad \quad \quad $\displaystyle R(x)=R_{m,\epsilon }(x),\ w_{1}(x):=R(x)^{r\delta
 },\ w_{2}(x):=R(x)^{r\gamma },\ w_{3}(x):=R(x)^{r\beta },$\par 
with the notation in Definition~\ref{pBB29}. Then, for any $\alpha
 >0,\ r\geq 2,$ we have:\par 
\quad $\displaystyle \forall \omega \in L^{r}(\lbrack 0,T+\alpha \rbrack
 ,L^{r}_{G}(M,w_{3}))\cap L^{r}(\lbrack 0,T+\alpha \rbrack ,L^{2}_{G}(M)),\
 \exists u\in L^{r}(\lbrack 0,T\rbrack ,W^{m,r}_{G}(M))::Du=\omega ,$\par 
with\par 
\quad \quad \quad $\displaystyle {\left\Vert{\partial _{t}u}\right\Vert}_{L^{r}(\lbrack
 0,T\rbrack ,L^{r}_{G}(M,w_{1}))}+{\left\Vert{u}\right\Vert}_{L^{r}(\lbrack
 0,T\rbrack ,W^{m,r}_{G}(M,w_{2}))}\leq $\par 
\quad \quad \quad \quad \quad $\displaystyle \leq c_{1}{\left\Vert{\omega }\right\Vert}_{L^{r}(\lbrack
 0,T+\alpha \rbrack ,L^{r}_{G}(M,w_{3}))}+c_{2}{\left\Vert{\omega
 }\right\Vert}_{L^{r}(\lbrack 0,T+\alpha \rbrack ,L_{G}^{2}(M))},$\par 
In the case of functions instead of sections of $G$ we have the
 same estimates but with $R(x)=R_{m-1,\epsilon }(x)$ and the weights:
\end{cor}
\quad \quad \quad $\displaystyle w_{1}(x):=R(x)^{r\delta '},\ w_{2}(x):=R(x)^{r\gamma
 '},\ w_{3}(x):=R(x)^{r\beta '}.$\ \par 
\quad Proof.\ \par 
By the threshold hypothesis (THL2), because $\displaystyle \omega
 \in L^{r}(\lbrack 0,T+\alpha \rbrack ,L^{2}_{G}(M))$ there is
 a $\displaystyle u\in L^{r}(\lbrack 0,T\rbrack ,L^{2}_{G}(M))$
 such that $\displaystyle Du=\omega $ and:\ \par 
\quad \quad \quad \begin{equation} {\left\Vert{u}\right\Vert}_{L^{r}(\lbrack 0,T+\alpha
 \rbrack ,L_{G}^{2}(M))}\lesssim {\left\Vert{\omega }\right\Vert}_{L^{r}(\lbrack
 0,T+\alpha \rbrack ,L_{G}^{2}(M))}.\label{pBB18}\end{equation}\ \par 
Hence, using Theorem~\ref{pBB16}, we get that the \emph{same}
 $u$ verifies:\ \par 
\quad \quad \quad $\displaystyle {\left\Vert{\partial _{t}u}\right\Vert}_{L^{r}(\lbrack
 0,T\rbrack ,L_{G}^{r}(M,w_{1}))}+{\left\Vert{u}\right\Vert}_{L^{r}(\lbrack
 0,T\rbrack ,W_{G}^{m,r}(M,w_{2}))}\leq $\ \par 
\quad \quad \quad \quad \quad $\displaystyle \leq c_{1}(k){\left\Vert{Du}\right\Vert}_{L^{r}(\lbrack
 0,T+\alpha \rbrack ,L_{G}^{r}(M,w_{3}))}+c_{2}(k){\left\Vert{u}\right\Vert}_{L^{r}(\lbrack
 0,T+\alpha \rbrack ,L_{G}^{2}(M))}.$\ \par 
So replacing $Du$ by $\omega $ and using~(\ref{pBB18}) we get\ \par 
\quad \quad \quad $\displaystyle {\left\Vert{\partial _{t}u}\right\Vert}_{L^{r}(\lbrack
 0,T\rbrack ,L_{G}^{r}(M,w_{1}))}+{\left\Vert{u}\right\Vert}_{L^{r}(\lbrack
 0,T\rbrack ,W_{G}^{m,r}(M,w_{2}))}\leq $\ \par 
\quad \quad \quad \quad \quad $\displaystyle \leq c_{1}{\left\Vert{\omega }\right\Vert}_{L^{r}(\lbrack
 0,T+\alpha \rbrack ,L_{G}^{r}(M,w_{3}))}+c_{2}{\left\Vert{\omega
 }\right\Vert}_{L^{r}(\lbrack 0,T+\alpha \rbrack ,L_{G}^{2}(M))}.$\ \par 
\quad The results for functions instead of sections of $G,$ follow
 the same lines and we have the same estimates but with $R(x)=R_{m-1,\epsilon
 }(x)$ and the weights:\ \par 
\quad \quad \quad $\displaystyle w_{1}(x):=R(x)^{r\delta '},\ w_{2}(x):=R(x)^{r\gamma
 '},\ w_{3}(x):=R(x)^{r\beta '}.$\ \par 
\quad The proof is complete. $\blacksquare $\ \par 
\ \par 
\quad If we are more interested in $L^{r}-L^{s}$ estimates, we can
 use the Sobolev embedding Theorem with weights~\cite{SobAmar19},
 valid here, which gives:\ \par 

\begin{thm}
Let $(M,g)$ be a complete riemannian manifold. Let $\displaystyle
 w(x):=R(x)^{\alpha }$ and $\displaystyle w':=R(x)^{\nu }$ with
  $\displaystyle \nu :=s(2+\alpha /r).$ Then $\displaystyle W_{G}^{m,r}(M,w)$
 is embedded in $\displaystyle W_{G}^{k,s}(M,w'),$ with $\displaystyle
 \frac{1}{s}=\frac{1}{r}-\frac{(m-k)}{n}>0$ and:\par 
\quad \quad \quad $\displaystyle \forall u\in W_{G}^{m,r}(M,w),\ {\left\Vert{u}\right\Vert}_{W_{G}^{k,s}(M,w')}\leq
 C{\left\Vert{u}\right\Vert}_{W_{G}^{m,r}(M,w)}.$
\end{thm}
So, with $\frac{1}{s}=\frac{1}{r}-\frac{(m-k)}{n}>0$ and $k=0$
 we get $\frac{1}{s}=\frac{1}{r}-\frac{m}{n}>0$ i.e. $s=\frac{nr}{n-rm}.$\ \par 
so $\displaystyle w(x):=R(x)^{b}\Rightarrow w':=R(x)^{\nu }$ with:\ \par 
\quad \quad \quad $\displaystyle \frac{\nu }{s}=2+\frac{b}{r}\Rightarrow \frac{\nu
 }{r}-m\frac{\nu }{n}=2+\frac{b}{r}\Rightarrow \nu (\frac{1}{r}-\frac{m}{n})=\frac{2r+b}{r}\Rightarrow
 \nu (n-m)=n(2r+b)$\ \par 
and finally $\displaystyle \nu =\frac{n(2r+b)}{n-m}.$ So we get:\ \par 

\begin{cor}
Let $M$ be a complete riemannian manifold of class ${\mathcal{C}}^{m}$
 without boundary. Moreover suppose we have (THL2). Then, on
 an adapted vector bundle $G,$ with $\displaystyle \nu =\frac{n(2r+\gamma
 )}{n-m}$ and $w_{4}(x):=R(x)^{r\nu }$ and also $s=\frac{nr}{n-rm}$:\par 
\quad $\displaystyle \forall \omega \in L^{r}(\lbrack 0,T+\alpha \rbrack
 ,L^{r}_{G}(M,w_{3}))\cap L^{r}(\lbrack 0,T+\alpha \rbrack ,L^{2}_{G}(M)),\
 \exists u\in L^{r}(\lbrack 0,T\rbrack ,L^{s}_{G}(M))::Du=\omega ,$\par 
such that:\par 
\quad \quad \quad $\displaystyle {\left\Vert{\partial _{t}u}\right\Vert}_{L^{r}(\lbrack
 0,T\rbrack ,L_{G}^{r}(M,w_{1}))}+{\left\Vert{u}\right\Vert}_{L^{r}(\lbrack
 0,T\rbrack ,L_{G}^{s}(M,w_{4}))}\leq $\par 
\quad \quad \quad \quad \quad $\displaystyle \leq c_{1}{\left\Vert{\omega }\right\Vert}_{L^{r}(\lbrack
 0,T+\alpha \rbrack ,L_{G}^{r}(M,w_{3}))}+c_{2}{\left\Vert{\omega
 }\right\Vert}_{L^{r}(\lbrack 0,T+\alpha \rbrack ,L_{G}^{2}(M))},$\par 
with $w_{1}(x):=R(x)^{r\delta },\ w_{3}(x):=R(x)^{r\beta },\
 w_{4}(x):=R(x)^{r\nu }.$\par 
\quad In the case of functions instead of sections of $G$ we have the
 same estimates but with $R(x)=R_{m-1,\epsilon }(x)$ and\par 
\quad \quad \quad $\displaystyle w_{1}(x):=R(x)^{r\delta '},\ w_{3}(x):=R(x)^{r\beta
 '},\ w_{4}(x):=R(x)^{r\nu }.$
\end{cor}

\subsection{The heat equation.}
\quad We shall consider the heat equation, $Du:=\partial _{t}u+\Delta
 u=\omega ,$ with $\Delta :=dd^{*}+d^{*}d$ the Hodge laplacian.
 Here we change the sign to use the standard notation with $\Delta
 $ essentially positive.\ \par 
\quad In this section we shall only consider the vector bundle of $p$-forms
 on the riemannian manifold $M.$ We denote $L^{r}_{p}(M)$ the
 space of $p$-forms in $L^{r}(M).$ The same for $\displaystyle
 W^{k,r}_{p}(M),$ the Sobolev spaces of $p$-forms on $M.$\ \par 
\quad We get that $\Delta ,$ the Hodge laplacian, is a $(C,\theta )$-elliptic
 operator on the $p$-forms in a complete riemannian manifold,
 for any $\theta >0,$ because its spectrum is contained in ${\mathbb{R}}_{+}$
 and\ \par 
\quad \quad \quad $\displaystyle \forall x\in M,\ \forall \xi _{x}\in T_{x}^{*}(M),\
 \left\vert{\xi _{x}}\right\vert =1,\ {\left\Vert{\Delta (x,\xi
 _{x})^{-1}}\right\Vert}\leq C.$\ \par 
\quad By Theorem~\ref{pBB19} we also have that the (THL2) hypothesis
 is true in this case, so we can apply Corollary~\ref{pBB20} to get:\ \par 

\begin{thm}
Let $M$ be a connected complete $n$-dimensional ${\mathcal{C}}^{2}$
 riemannian manifold without boundary. Let  $\displaystyle Du:=\partial
 _{t}u+\Delta u$ be the heat operator acting on the bundle $\displaystyle
 \Lambda ^{p}(M)$ of $p$-forms on $M.$ Let:\par 
\quad \quad \quad $\displaystyle R(x)=R_{2,\epsilon }(x),\ w_{1}(x):=R(x)^{r\delta
 },\ w_{2}(x):=R(x)^{r\gamma },\ w_{3}(x):=R(x)^{r\beta },$\par 
with the notation in Definition~\ref{pBB29} with $m=2.$ Then,
 for any $\alpha >0,\ r\geq 2,$ we have:\par 
\quad $\displaystyle \forall \omega \in L^{r}(\lbrack 0,T+\alpha \rbrack
 ,L^{r}_{p}(M,w_{3}))\cap L^{r}(\lbrack 0,T+\alpha \rbrack ,L^{2}_{p}(M)),\
 \exists u\in L^{r}(\lbrack 0,T\rbrack ,W^{2,r}_{p}(M))::Du=\omega ,$\par 
with\par 
\quad \quad \quad $\displaystyle {\left\Vert{\partial _{t}u}\right\Vert}_{L^{r}(\lbrack
 0,T\rbrack ,L^{r}_{p}(M,w_{1}))}+{\left\Vert{u}\right\Vert}_{L^{r}(\lbrack
 0,T\rbrack ,W^{2,r}_{p}(M,w_{2}))}\leq $\par 
\quad \quad \quad \quad \quad $\displaystyle \leq c_{1}{\left\Vert{\omega }\right\Vert}_{L^{r}(\lbrack
 0,T+\alpha \rbrack ,L^{r}_{p}(M,w_{3}))}+c_{2}{\left\Vert{\omega
 }\right\Vert}_{L^{r}(\lbrack 0,T+\alpha \rbrack ,L_{p}^{2}(M))},$\par 
In the case of functions i.e. $p=0,$ we have the same estimates
 but with $R(x)=R_{1,\epsilon }(x)$ and\par 
\quad \quad \quad $\displaystyle w_{1}(x):=R(x)^{r\delta '},\ w_{2}(x):=R(x)^{r\gamma
 '},\ w_{3}(x):=R(x)^{r\beta '}.$
\end{thm}

\section{Classical estimates.}
\quad We shall give some examples where we have classical estimates
 using that $\displaystyle \forall x\in M,\ R_{\epsilon }(x)\geq
 \delta ,$ via~\cite[Corollary, p. 7] {HebeyHerzlich97}  (see
 also Theorem 1.3 in the book by Hebey~\cite{Hebey96}):\ \par 

\begin{cor}
~\label{pBB17}Let $(M,g)$ be a complete riemannian manifold.
 Let $m\geq 1$; if we have the injectivity radius $\displaystyle
 r_{inj}(x)\geq i>0$ and  $\displaystyle \forall j\leq m-1,\
 \left\vert{\nabla ^{j}Rc_{(M,g)}(x)}\right\vert \leq c$ for
 all $\displaystyle x\in M,$ then there exists a constant $\delta
 >,0,$ depending only on $\displaystyle n,\epsilon ,i,m$ and
  $\displaystyle c,$ such that: $\displaystyle \forall x\in M,\
 R_{m,\epsilon }(x)\geq \delta .$
\end{cor}
\quad Proof.\ \par 
The Theorem of Hebey and Herzlich gives that, under these hypotheses,
 for any $\alpha \in (0,1)$ there exists a constant $\delta >0,$
 depending only on $\displaystyle n,\epsilon ,i,m,\alpha $ and
  $\displaystyle c,$ such that:\ \par 
\quad \quad \quad $\displaystyle \forall x\in M,\ r_{H}(1+\epsilon ,m,\alpha )(x)\geq
 \delta >0.$\ \par 
So even taking our definition with a harmonic coordinates patch,
 we have that:\ \par 
\quad \quad \quad $\displaystyle R_{m,\epsilon }(x)\geq r_{H}(1+\epsilon ,m,\alpha )(x)$\ \par 
so a fortiori when we take the sup for $R_{m,\epsilon }(x)$ on
 \emph{any} smooth coordinates patch, not necessarily harmonic patch.\ \par 
\quad The proof is complete. $\blacksquare $\ \par 

\subsection{Bounded geometry.}

\begin{defin}
A riemannian manifold $M$ has $k$-order {\bf bounded geometry} if:\par 
\quad $\bullet $  the injectivity radius $r_{inj}(x)$ at $x\in M$ is
 bounded below by some constant $\delta >0$ for any $\displaystyle x\in M$\par 
\quad $\bullet $ and if for $0\leq j\leq k,$ the covariant derivatives
 $\nabla ^{j}R$ of the curvature tensor are bounded in $L^{\infty }(M)$ norm.
\end{defin}
We shall weakened this definition to suit our purpose.\ \par 

\begin{defin}
A riemannian manifold $M$ has $k$-order  {\bf weak bounded geometry} if:\par 
\quad $\bullet $  the injectivity radius $r_{inj}(x)$ at $x\in M$ is
 bounded below by some constant $\delta >0$ for any $\displaystyle x\in M$\par 
\quad $\bullet $ and if for $0\leq j\leq k,$ the covariant derivatives
 $\nabla ^{j}Rc$ of the Ricci curvature tensor are bounded in
 $L^{\infty }(M)$ norm.
\end{defin}
\quad Using this notion, we get our main Theorem~\ref{pBB16} without weights:\ \par 

\begin{thm}
~\label{pBL41}Suppose that $A$ is a $(C,\theta )$-elliptic operator
 of order $m$ acting on sections of the adapted vector bundle
 $\displaystyle G:=(H,\pi ,M)$ in the complete riemannian manifold
 $(M,g),$ with $\theta <\pi /2,$ and consider the parabolic equation
 $\displaystyle Du=\partial _{t}u-Au$ also acting on sections
 of $G.$ Suppose moreover that $(M,g)$ has $(m-1)$ order weak
 bounded geometry and (THL2) is true. Then\par 
\quad $\displaystyle \forall \omega \in L^{r}(\lbrack 0,T+\alpha \rbrack
 ,L^{r}_{G}(M))\cap L^{r}(\lbrack 0,T+\alpha \rbrack ,L^{2}_{G}(M)),\
 \exists u\in L^{r}(\lbrack 0,T\rbrack ,W^{m,r}_{G}(M))::Du=\omega ,$\par 
with:\par 
\quad $\displaystyle {\left\Vert{\partial _{t}u}\right\Vert}_{L^{r}(\lbrack
 0,T\rbrack ,L_{G}^{r}(M))}+{\left\Vert{u}\right\Vert}_{L^{r}(\lbrack
 0,T\rbrack ,W_{G}^{m,r}(M))}\leq c_{1}{\left\Vert{\omega }\right\Vert}_{L^{r}(\lbrack
 0,T+\alpha \rbrack ,L_{G}^{r}(M))}+c_{2}{\left\Vert{\omega }\right\Vert}_{L^{r}(\lbrack
 0,T+\alpha \rbrack ,L_{G}^{2}(M))}.$\par 
In the case of functions instead of sections of $G$ we have the
 same estimates just supposing that $(M,g)$ has $(m-2)$ order
 weak bounded geometry.
\end{thm}

\subsubsection{Examples of manifolds of bounded geometry.}
\quad $\bullet $  Euclidean space with the standard metric has bounded
 geometry.
\ \par 
\ \par 
\quad $\bullet $ A smooth, compact Riemannian manifold $M$ has bounded
 geometry as well;
 both the injectivity radius and the curvature
 including derivatives are continuous
 functions, so these attain
 their finite minima and maxima, respectively on $M.$ If $M\in
 {\mathcal{C}}^{m+2},$ then it has bounded geometry of order $m.$\ \par 
\ \par 
\quad $\bullet $ Non compact, smooth Riemannian manifolds that possess
 a transitive group
 of isomorphisms (such as the hyperbolic
 spaces ${\mathbb{H}}^{n}$) have $m$-order bounded geometry since
 the
 finite injectivity radius and curvature estimates at any
 single point translate to
 a uniform estimate for all points
 under isomorphisms.\ \par 
\ \par 
\quad Of course these examples have a fortiori weak bounded geometry.\ \par 

\subsection{Hyperbolic manifolds.~\label{p37}}
\quad These are manifolds such that the sectional curvature $K_{M}$
 is constantly $-1.$ For them we have first that the Ricci curvature
 is bounded.\ \par 

\begin{lem}
~\label{BG0} Let $(M,g)$ be a complete Riemannian manifold such
 that $H\leq K_{M}\leq K$ for constants
 $H,K\in {\mathbb{R}}.$
 Then we have that ${\left\Vert{Rc}\right\Vert}_{\infty }\leq
 \max (\left\vert{H}\right\vert ,\left\vert{K}\right\vert ).$
\end{lem}
\quad This lemma is so well known than we can omit its proof. $\blacksquare $\ \par 
\ \par 
\quad To get that the injectivity radius $r_{inj}(x)$ is bounded below
 we shall use a Theorem by Cheeger, Gromov
 and Taylor~\cite{CheegerGroTay82}:\
 \par 

\begin{thm}
~\label{pL33}Let $(M,g)$ be a complete Riemannian manifold such
 that $K_{M}\leq K$ for constants
 $K\in {\mathbb{R}}.$ Let $0<r<\frac{\pi
 }{4{\sqrt{K}}}$ if $K>0$ and $r\in (0,\infty )$ if $K\leq 0.$
 Then the injectivity radius $r_{inj}(x)$ at $x$ satisfies\par 
\quad \quad \quad $\displaystyle r_{inj}(x)\geq r\frac{\mathrm{V}\mathrm{o}\mathrm{l}(B_{M}(x,r))}{\mathrm{V}\mathrm{o}\mathrm{l}(B_{M}(x,r))+\mathrm{V}\mathrm{o}\mathrm{l}(B_{T_{x}M}(0,2r))},$\par
 
where $\displaystyle B_{T_{x}M}(0,2r))$ denotes the volume of
 the ball of radius $2r$ in $T_{x}M,$ where both the
 volume
 and the distance function are defined using the metric $g^{*}:=\exp
 _{p}^{*}g$ i.e. the pull-back of
 the metric $g$ to $\displaystyle
 T_{x}M$ via the exponential map.
\end{thm}
\quad This Theorem leads to the definition:\ \par 

\begin{defin}
Let $(M,g)$ be a Riemannian manifold. We shall say that it has
 the {\bf lifted doubling property} if we have:\par 
\quad $\displaystyle (LDP)\ \ \ \ \ \ \ \ \ \exists \alpha ,\beta >0::\forall
 x\in M,\ \exists r\geq \beta ,\ \mathrm{V}\mathrm{o}\mathrm{l}(B_{T_{x}M}(0,2r))\leq
 \alpha \mathrm{V}\mathrm{o}\mathrm{l}(B_{M}(x,r)),$\par 
where $\displaystyle B_{T_{x}M}(0,2r))$ denotes the volume of
 the ball of radius $2r$ in $T_{x}M,$ and both the
 volume and
 the distance function are defined on $T_{x}M$ using the metric
 $g^{*}:=\exp _{p}^{*}g$ i.e. the pull-back of
 the metric $g$
 to $\displaystyle T_{x}M$ via the exponential map.
\end{defin}
Hence we get:\ \par 

\begin{cor}
~\label{pL44}Let $(M,g)$ be a complete Riemannian manifold such
 that $K_{M}\leq K$ for a constant
 $K\in {\mathbb{R}}.$ For
 instance an hyperbolic manifold.  Suppose moreover that $(M,g),$
 has the lifted doubling property.\par 
Then $\forall x\in M,\ r_{inj}(x)\geq \frac{\beta }{1+\alpha }.$
\end{cor}
\quad Proof.\ \par 
By the (LDP) we get, for a $r\geq \beta ,$\ \par 
\quad \quad \quad $\displaystyle \mathrm{V}\mathrm{o}\mathrm{l}(B_{T_{x}M}(0,2r))\leq
 \alpha \mathrm{V}\mathrm{o}\mathrm{l}(B_{M}(x,r)).$\ \par 
We apply Theorem~\ref{pL33} of Cheeger, Gromov
 and Taylor to get\ \par 
\quad \quad \quad $\displaystyle r_{inj}(x)\geq r\frac{\mathrm{V}\mathrm{o}\mathrm{l}(B_{M}(x,r))}{\mathrm{V}\mathrm{o}\mathrm{l}(B_{M}(x,r))+\mathrm{V}\mathrm{o}\mathrm{l}(B_{T_{x}M}(0,2r))}.$\
 \par 
So\ \par 
\quad \quad \quad $\displaystyle \frac{\mathrm{V}\mathrm{o}\mathrm{l}(B_{M}(x,r))}{\mathrm{V}\mathrm{o}\mathrm{l}(B_{M}(x,r))+\mathrm{V}\mathrm{o}\mathrm{l}(B_{T_{x}M}(0,2r))}\geq
 \frac{1}{1+\alpha }$\ \par 
hence, because $\displaystyle r\geq \beta ,$ we get the result.
 $\blacksquare $\ \par 
\ \par 
\quad As an example of application we get\ \par 

\begin{prop}
~\label{pBL40}Let $(M,g)$ be a complete Riemannian manifold such
 that $H\leq K_{M}\leq K$ for constants
 $H,K\in {\mathbb{R}},$
 where $\displaystyle K_{M}$ is the sectional curvature of $M.$
 Suppose moreover that  $(M,g)$ has the lifted doubling property
 and that, for $0\leq j\leq k,$ the covariant derivatives $\nabla
 ^{j}Rc$ of the Ricci curvature tensor are bounded in $L^{\infty
 }(M)$ norm. Then $(M,g)$ has weak bounded geometry of order $k.$
\end{prop}

\begin{thm}
Let $(M,g)$ be a complete Riemannian manifold such that $H\leq
 K_{M}\leq K$ for constants
 $H,K\in {\mathbb{R}},$ where $\displaystyle
 K_{M}$ is the sectional curvature of $M.$ Suppose moreover that
  $(M,g)$ has the lifted doubling property. Suppose that $A$
 is a $(C,\theta )$-elliptic operator of order $m$ acting on
 sections of the adapted vector bundle $\displaystyle G:=(H,\pi
 ,M)$ in $(M,g),$ with $\theta <\pi /2,$ and consider the parabolic
 equation $\displaystyle Du=\partial _{t}u-Au$ also acting on
 sections of $G.$ Moreover suppose we have (THL2). Provided that,
 for $0\leq j\leq m-1,$ the covariant derivatives $\nabla ^{j}Rc$
 of the Ricci curvature tensor are bounded in $L^{\infty }(M)$ norm:\par 
\quad $\displaystyle \forall \omega \in L^{r}(\lbrack 0,T+\alpha \rbrack
 ,L^{r}_{G}(M))\cap L^{r}(\lbrack 0,T+\alpha \rbrack ,L^{2}_{G}(M)),\
 \exists u\in L^{r}(\lbrack 0,T\rbrack ,W^{m,r}_{G}(M))::Du=\omega ,$\par 
with:\par 
\quad $\displaystyle {\left\Vert{\partial _{t}u}\right\Vert}_{L^{r}(\lbrack
 0,T\rbrack ,L_{G}^{r}(M))}+{\left\Vert{u}\right\Vert}_{L^{r}(\lbrack
 0,T\rbrack ,W_{G}^{m,r}(M))}\leq c_{1}{\left\Vert{\omega }\right\Vert}_{L^{r}(\lbrack
 0,T+\alpha \rbrack ,L_{G}^{r}(M))}+c_{2}{\left\Vert{\omega }\right\Vert}_{L^{r}(\lbrack
 0,T+\alpha \rbrack ,L_{G}^{2}(M))}.$\par 
In the case of functions instead of sections of $G$ we have the
 same estimates, just supposing that for $0\leq j\leq m-2,$ the
 covariant derivatives $\nabla ^{j}Rc$ of the Ricci curvature
 tensor are bounded in $L^{\infty }(M)$ norm.
\end{thm}
\quad Proof.\ \par 
By the Proposition~\ref{pBL40}, we have that $(M,g)$ has weak
 bounded geometry of order $k.$ So we can apply Theorem~\ref{pBL41}.
 $\blacksquare $\ \par 
\ \par 
\quad And in the case of the heat equation:\ \par 

\begin{cor}
Let $(M,g)$ be a complete Riemannian manifold such that $H\leq
 K_{M}\leq K$ for constants
 $H,K\in {\mathbb{R}},$ where $\displaystyle
 K_{M}$ is the sectional curvature of $M.$ Suppose moreover that
  $(M,g)$ has the lifted doubling property. Then $\displaystyle
 \exists \delta >0,\ \forall x\in M,\ R_{1,\epsilon }(x)\geq
 \delta .$ This implies that we get "classical solutions" for
 the heat equation for functions in this case. I.e.\par 
\quad $\displaystyle \forall \omega \in L^{r}(\lbrack 0,T+\alpha \rbrack
 ,L^{r}(M))\cap L^{r}(\lbrack 0,T+\alpha \rbrack ,L^{2}(M)),\
 \exists u\in L^{r}(\lbrack 0,T\rbrack ,W^{2,r}(M))::Du=\omega ,$\par 
with:\par 
\quad $\displaystyle {\left\Vert{\partial _{t}u}\right\Vert}_{L^{r}(\lbrack
 0,T\rbrack ,L^{r}(M))}+{\left\Vert{u}\right\Vert}_{L^{r}(\lbrack
 0,T\rbrack ,W^{2,r}(M))}\leq c_{1}{\left\Vert{\omega }\right\Vert}_{L^{r}(\lbrack
 0,T+\alpha \rbrack ,L^{r}(M))}+c_{2}{\left\Vert{\omega }\right\Vert}_{L^{r}(\lbrack
 0,T+\alpha \rbrack ,L^{2}(M))}.$\par 
\quad To get that $\displaystyle \exists \delta >0,\ \forall x\in M,\
 R_{2,\epsilon }(x)\geq \delta ,$ we need to ask that moreover
 the covariant derivatives $\nabla Rc$ of the Ricci curvature
 tensor are bounded in $L^{\infty }(M)$ norm. This is the case
 in particular if $(M,g)$ is hyperbolic. This implies that we
 get "classical solutions" for the heat equation for $p$-forms
 in this case. I.e.\par 
\quad $\displaystyle \forall \omega \in L^{r}(\lbrack 0,T+\alpha \rbrack
 ,L^{r}_{p}(M))\cap L^{r}(\lbrack 0,T+\alpha \rbrack ,L^{2}_{p}(M)),\
 \exists u\in L^{r}(\lbrack 0,T\rbrack ,W^{2,r}_{p}(M))::Du=\omega ,$\par 
with:\par 
\quad $\displaystyle {\left\Vert{\partial _{t}u}\right\Vert}_{L^{r}(\lbrack
 0,T\rbrack ,L^{r}_{p}(M))}+{\left\Vert{u}\right\Vert}_{L^{r}(\lbrack
 0,T\rbrack ,W^{2,r}_{p}(M))}\leq c_{1}{\left\Vert{\omega }\right\Vert}_{L^{r}(\lbrack
 0,T+\alpha \rbrack ,L^{r}_{p}(M))}+c_{2}{\left\Vert{\omega }\right\Vert}_{L^{r}(\lbrack
 0,T+\alpha \rbrack ,L^{2}_{p}(M))}.$
\end{cor}
\quad Proof.\ \par 
By Lemma~\ref{BG0} we get that $\displaystyle {\left\Vert{Rc}\right\Vert}_{\infty
 }<\infty .$ Then we apply Corollary~\ref{pL44}. For forms we
 have to use the extra hypothesis on the covariant derivatives.
 $\blacksquare $\ \par 

\begin{rem}
In the case the hyperbolic manifold $(M,g)$ is simply connected,
 then by the Hadamard Theorem~\cite[Theorem 3.1, p. 149]{DoCarmo93},
 we get that the injectivity radius is $\infty ,$ so we have
 also the classical estimates in this case.
\end{rem}
\ \par 
\ \par 

\bibliographystyle{/usr/local/texlive/2017/texmf-dist/bibtex/bst/base/apalike}

\end{document}